\documentclass[12pt]{amsart}

\usepackage{a4wide}

\usepackage{amssymb,amsmath,stmaryrd,xy}

\usepackage{enumerate}

\newtheorem{theorem}{Theorem}[section]
\newtheorem{proposition}[theorem]{Proposition}
\newtheorem{lemma}[theorem]{Lemma}
\newtheorem{corollary}[theorem]{Corollary}

\newtheorem*{maintheorem}{Main Theorem}

\theoremstyle{remark}
\newtheorem{remark}[theorem]{\it Remark}

\newcommand{\C}{\ensuremath{\mathbb{C}}}
\newcommand{\R}{\ensuremath{\mathbb{R}}}

\newcommand{\g}[1]{\ensuremath{\mathfrak{#1}}}

\newcommand{\II}{\ensuremath{I\!I}}
\newcommand{\sech}{\ensuremath{\mathop{\rm sech}\nolimits}}
\newcommand{\csch}{\ensuremath{\mathop{\rm csch}\nolimits}}

\begin{document}
\title[Non-Hopf hypersurfaces with constant principal curvatures]
{Non-Hopf real hypersurfaces with constant\\ principal curvatures in
complex space forms}

\author[J. C. D\'\i{}az-Ramos]{Jos\'{e} Carlos D\'\i{}az-Ramos}
\address{Department of Geometry and Topology,
University of Santiago de Compostela, Spain.}
\email{josecarlos.diaz@usc.es}

\author[M. Dom\'{\i}nguez-V\'{a}zquez]{Miguel Dom\'{\i}nguez-V\'{a}zquez}
\address{Department of Geometry and Topology,
University of Santiago de Compostela, Spain.}
\email{miguel.dominguez@usc.es}

\thanks{The first author has been supported by a Marie-Curie
European Reintegration Grant (PERG04-GA-2008-239162). The
second author has been supported by the FPU programme of the
Spanish Government. Both authors have been supported by project
MTM2006-01432 (Spain).}

\begin{abstract}
We classify real hypersurfaces in complex space forms with
constant principal curvatures and whose Hopf vector field has two
nontrivial projections onto the principal curvature spaces.

In complex projective spaces
such real hypersurfaces do not exist. In complex hyperbolic
spaces these are holomorphically congruent to open parts of
tubes around the ruled minimal submanifolds with totally real
normal bundle introduced by Berndt and Br\"{u}ck.
In particular, they are open parts of homogenous ones.
\end{abstract}


\phantom{Draft version, \today}

\maketitle

%\begin{center}{Draft version, \today}\end{center}

%%%%%%%%%%%%%%%%%%%%%%%%%% Section %%%%%%%%%%%%%%%%%%%%%%%%%%%%%%%
\section{Introduction}

A homogeneous submanifold of a Riemannian manifold is an orbit
of the action of a closed subgroup of the isometry group of the
ambient manifold. One of the aims of submanifold geometry is to
classify homogeneous submanifolds of a given
manifold and to characterize them in terms of geometric data. Of particular interest are homogeneous
hypersurfaces, which arise as principal orbits of cohomogeneity
one actions. Obviously,
homogeneous hypersurfaces have constant principal curvatures, that is, the eigenvalues of their shape operator are constant.
It is an outstanding problem to determine under
which conditions  hypersurfaces with constant principal
curvatures are open parts of homogeneous ones.

In spaces of constant curvature, a hypersurface has constant
principal curvatures if and only if it is isoparametric. The classification of isoparametric hypersurfaces was achieved by Segre
\cite{Se38} in Euclidean
spaces and by Cartan \cite{Ca38} in real hyperbolic spaces. They all are open parts of homogeneous
ones. The situation is more involved in
spheres. Cartan classified hypersurfaces with $g\in\{1,2,3\}$ constant
principal curvatures in spheres. Subsequently, Hsiang and Lawson \cite{HL71} classified homogeneous hypersurfaces in spheres; they have $g\in\{1,2,3,4,6\}$ principal curvatures. Then, M\"{u}nzner \cite{Mu80} showed that $g\in\{1,2,3,4,6\}$
for isoparametric hypersurfaces in general. Surprisingly, for
$g=4$ there are isoparametric hypersurfaces that are not homogeneous \cite{FKM81}. Recently, Cecil, Chi and Jensen \cite{CCJ07}, and Immervoll \cite{I08} showed that, with a few
possible exceptions, hypersurfaces with $g=4$ constant
principal curvatures are among the known homogeneous and
inhomogeneous examples. Some progress has been made for $g=6$
by Abresch \cite{Ab83} and Dorfmeister and Neher \cite{DN85},
but the problem remains open in full generality. See \cite{Th00} for a
survey.

The problem is even more difficult in complex space forms. See
\cite{NR97} for a survey on this and related topics. By $c\neq 0$ we
denote the constant holomorphic sectional curvature of a
complex space form; thus, if $c>0$ (resp. $c<0$) we have a complex projective (resp. hyperbolic)
space $\C P^n(c)$ (resp. $\C H^n(c)$). We denote by $J$
its K\"{a}hler structure. Let $M$ be a real hypersurface of a complex
space form and $\xi$ a (local) unit normal vector field. Then,
$J\xi$ is tangent to $M$ and is called the Hopf vector field of $M$. The hypersurface $M$ is said to be Hopf if $J\xi$ is
a principal curvature vector field. The motivation for our work is
to address the classification of real hypersurfaces with constant principal curvatures in complex space forms. We briefly summarize the current
state of the problem.

Assume $M$ is a real hypersurface of a complex space form with $g$ distinct constant principal curvatures. For $p\in M$ denote by $h(p)$ the number of nontrivial projections of $J\xi_p$ onto the principal curvature spaces of
$M$. Clearly, this function is integer-valued and $M$ is Hopf
if and only if $h=1$. The classification of homogeneous real hypersurfaces in complex projective spaces $\C P^n(c)$ was
derived by Takagi \cite{Ta73}. It follows from this
classification that $g\in\{2,3,5\}$. A remarkable feature of
homogeneous real hypersurfaces in $\C P^n(c)$ is that they are
Hopf. Subsequently, Takagi classified real hypersurfaces with
$g\in\{2,3\}$ constant principal curvatures \cite{Ta75a},
\cite{Ta75b} (\cite{W83} for $n=2$, $g=3$). It follows from his work that they all are Hopf
and open parts of homogeneous ones. Kimura \cite{Ki86}
classified Hopf real hypersurfaces with constant principal
curvatures in $\C P^n$ and showed that these are open parts of
homogeneous ones. No examples are known of real hypersurfaces
with constant principal curvatures in $\C P^n(c)$ with $h>1$.
Surprisingly, in $\C H^n(c)$
there are non-Hopf homogeneous real hypersurfaces. The first example was discovered by Lohnherr \cite{L99} and further examples were given by Berndt and Br\"{u}ck \cite{Be98},
\cite{BB01}. We refer to \S\ref{secExamples} for a brief
introduction and to \cite{BD09} for a
deeper study of their geometry. Berndt and Tamaru obtained in \cite{BT07}
the classification of cohomogeneity one actions on $\C H^n(c)$.
The number of principal curvatures of the homogeneous examples
is $g\in\{2,3,4,5\}$. Montiel \cite{Mo85} classified real
hypersurfaces with $g=2$ constant principal curvatures in $\C
H^n(c)$ ($n\geq 3$). Berndt and the first author solved the
case $g=3$ and $g=2$, $n=2$ \cite{BD05}, \cite{BD06}. It
follows from \cite{Mo85} that $h=1$ when $g=2$, and from
\cite{BD05} and \cite{BD06} we get $h\leq 2$ if $g=3$. Hopf real
hypersurfaces with constant principal curvatures in $\C H^n(c)$
were classified by Berndt \cite{Be89} and they all are open
parts of homogeneous ones.  To our knowledge,
\cite{BD05} and \cite{BD06} are the first classifications of
this kind involving non-Hopf real hypersurfaces. Nothing is
known about $h$ if $g\geq 4$.

Our aim in this paper is to carry out the next natural step
after Berndt and Kimura's classification of Hopf real
hypersurfaces with constant principal curvatures in $\C P^n(c)$
and $\C H^n(c)$. Thus, we classify
real hypersurfaces with constant principal curvatures whose Hopf vector field $J\xi$ has $h=2$
nontrivial projections onto the principal curvature spaces.

\begin{maintheorem}\label{thMain}
We have:
\begin{enumerate}[{\rm (a)}]
\item\label{thMainCPn} There are no real hypersurfaces with
    constant principal curvatures in $\C P^n(c)$, $n\geq
    2$, whose Hopf vector field  has $h=2$ nontrivial projections onto the principal curvature spaces.

\item\label{thMainCHn} Let $M$ be a connected real
    hypersurface in $\C H^n(c)$, $n\geq 2$, with constant
    principal curvatures and whose Hopf vector field has $h=2$ nontrivial projections onto the \text{principal}
    curvature spaces of $M$. Then, $M$ has
    $g\in\{3,4\}$ principal curvatures and is
    holomorphically congruent to an open part of:
    \begin{itemize}
    \item[(i)] a ruled minimal real
        hypersurface $W^{2n-1}\subset \C H^n(c)$ or one
        of the equidistant hypersurfaces to $W^{2n-1}$,
        or
    \item[(ii)] a tube around a ruled minimal
        Berndt-Br\"{u}ck submanifold $W^{2n-k}\subset\C H^n(c)$ with totally real normal bundle, for some
        $k\in\{2,\dots,n-1\}$.
    \end{itemize}
    In particular, $M$ is an open part of a homogeneous
    real hypersurface of $\C H^n(c)$.
\end{enumerate}
\end{maintheorem}

The ruled minimal submanifolds $W^{2n-k}\subset\C
H^n(c)$ are homogeneous and have totally real normal bundle of
rank $k\in\{1,\dots,n-1\}$. Actually, $W^{2n-1}$ was discovered by Lohnherr \cite{L99}. Then, Berndt studied the geometry of the equidistant hypersurfaces to $W^{2n-1}$ \cite{Be98}. This construction was generalized by Berndt and Br\"{u}ck in \cite{BB01}. Both $W^{2n-1}$ and any of its
equidistant hypersurfaces have $g=3$ principal curvatures. The
tubes around $W^{2n-k}$, $k\in\{2,\dots,n-1\}$ have $g=4$
principal curvatures if
$r\neq({1}/{\sqrt{-c}})\log(2+\sqrt{3})$ and $g=3$ principal
curvatures if $r=({1}/{\sqrt{-c}})\log(2+\sqrt{3})$. See
\cite{BD09} for a detailed description.

The proof is as follows. First we use the Gauss and Codazzi
equations to derive some algebraic properties of the eigenvalue
structure of the shape operator. The methods used for this are similar to those of \cite{BD05}, although a bit more general. We
would like to emphasize that whenever we use a method similar
to one in \cite{BD05} we explicitly point it out and skip the
details as much as possible. On the other hand, we focus
on the new techniques and results, especially on Subsection~\ref{secBoundg}. The most
crucial step of the proof is to show that the number $g$ of
constant principal curvatures satisfies $g\leq 4$. For this we
use a novel approach based on the study of some inequalities
satisfied by the principal curvatures. Using standard Jacobi
field theory one can deduce the geometry of the focal
submanifolds of these hypersurfaces and then the result follows
from a rigidity result in \cite{BD09}.

The paper is organized as follows. In Section~\ref{secPreliminaries} we introduce the basic elements of our paper. Subsection~\ref{secEquations} is devoted to present the equations of submanifold geometry that we will use in the rest of the paper. In \S\ref{secExamples} we briefly describe the ruled minimal Berndt-Br\"{u}ck submanifolds $W^{2n-k}$. We prove our Main Theorem in Section~\ref{secProof}. The proof is divided in several steps. Some vector fields and functions arise naturally in our proof (\S\ref{secNotation} and \S\ref{secA}). We get some of their properties in Subsection~\ref{secPrincipalSpaces}. In \S\ref{secBoundg} we show that the number $g$ of principal curvatures satisfies $g\in\{3,4\}$. We summarize all the eigenvalue structure in \S\ref{secEigenvalueStructure}. In Subsection~\ref{secJacobi} we use standard Jacobi field theory to finish the proof of the Main Theorem.

%%%%%%%%%%%%%%%%%%%%%%%%%% Section %%%%%%%%%%%%%%%%%%%%%%%%%%%%%%%
\section{Preliminaries}\label{secPreliminaries}

In this section we introduce the basic notation of this paper.
We write down the Gauss and Codazzi equations of a hypersurface
in a complex space form and derive some basic consequences. Then,
we briefly mention how the examples of the Main Theorem are
constructed.

%%%%%%%%%%%%%%%%%%%%%%%%%% Subsection %%%%%%%%%%%%%%%%%%%%%%%%%%%%%%%
\subsection{The equations of a hypersurface}\label{secEquations}

Let $\bar{M}(c)$ be a complex space form of constant
holomorphic sectional curvature $c\neq 0$ and complex dimension
$n$. If $c>0$ then $\bar{M}(c)$ is a complex
projective space $\C P^n(c)$ of constant holomorphic sectional curvature $c$. Analogously, if $c<0$ then $\bar{M}(c)$ is a complex
hyperbolic space $\C H^n(c)$. We denote by
$\langle\,\cdot\,,\,\cdot\,\rangle$ its inner product, by $J$
its K\"{a}hler structure, and by $\bar\nabla$ its Levi-Civita
connection. The curvature tensor is defined by
$\bar{R}(X,Y)=[\bar\nabla_X,\bar\nabla_Y]-\bar\nabla_{[X,Y]}$, so in this
case we have
\[
\bar{R}(X,Y)Z=\frac{c}{4}\left(\langle Y,Z\rangle X-\langle X,Z\rangle Y
+\langle JY,Z\rangle JX-\langle JX,Z\rangle JY
-2\langle JX,Y\rangle JZ\right).
\]

Let $M$ be a connected submanifold of $\bar{M}(c)$. We denote
by $\nabla$ and $R$ its Levi-Civita connection and its
curvature tensor respectively. By $TM$ and $\nu M$ we denote
the tangent and normal bundles of $M$. We use the symbol
$\Gamma(\cdot)$ to refer to the smooth sections of any vector bundle.
Let $X,Y,Z,W\in\Gamma(TM)$ and $\xi\in\Gamma(\nu M)$.

The second fundamental form $\II$ of $M$ is defined by the
Gauss formula as $\bar{\nabla}_X Y=\nabla_X Y+\II(X,Y)$. The
Weingarten formula is then written as $\bar\nabla_X\xi=-S_\xi
X+\nabla_X^\perp\xi$, where $S_\xi$ is the shape operator with
respect to $\xi$ and $\nabla^\perp$ is the induced normal
connection on $\nu M$. The second fundamental form and the
shape operator are related by $\langle S_\xi
X,Y\rangle=\langle\II(X,Y),\xi\rangle$.

Now let $M$ be a connected real hypersurface of $\bar{M}(c)$. The word `real' emphasizes the fact that the \emph{real} codimension is one. Fix
$\xi\in\Gamma(\nu M)$ a (local) unit normal vector field. We write $S$ instead
of $S_\xi$. The Gauss formula can be rewritten as
\[
\bar{\nabla}_X Y=\nabla_X Y+\langle SX,Y\rangle\xi,
\]
and hence, the Weingarten formula is $SX=-\bar\nabla_X\xi$.
Moreover, the Gauss and Codazzi equations for a hypersurface are
\begin{eqnarray*}
\langle\bar{R}(X,Y)Z,W\rangle
    &=&\langle{R}(X,Y)Z,W\rangle-\langle SY,Z\rangle\langle SX,W\rangle
    +\langle SX,Z\rangle\langle SY,W\rangle,\mbox{ and}\\
\langle\bar{R}(X,Y)Z,\xi\rangle &=&\langle(\nabla_X S)Y-(\nabla_Y S)X,Z\rangle.
\end{eqnarray*}

We assume from now on that $M$ has constant principal
curvatures, that is, the eigenvalues of the shape operator $S$ are constant. For each principal curvature $\lambda$ of $M$ we
denote by $T_\lambda$ the distribution on $M$ formed by the
principal curvature spaces of $\lambda$ along $M$.

The Codazzi equation implies (see \cite [Section 2]{BD05} for a proof)

\begin{lemma}\label{thCodazzi}
\begin{enumerate}[{\rm (i)}]
\item\label{thCodazzi1Space} Let $p\in M$. If the
    orthogonal projection of $J\xi_p$ onto $T_\alpha(p)$ is
    nonzero, then $T_\alpha(p)$ is a real subspace of
    $T_p\bar{M}(c)$, that is, $JT_\alpha(p)$ is orthogonal
    to $T_\alpha(p)$.

\item\label{thCodazzi2Spaces} Let $X,Y\in\Gamma(T_\alpha)$
    and $Z\in\Gamma(T_\beta)$ with $\alpha\neq\beta$. Then
\[
\langle\nabla_X Y,Z\rangle
=\frac{c}{4(\alpha-\beta)}\left(\langle JY,Z\rangle\langle X,J\xi\rangle
+\langle JX,Y\rangle\langle Z,J\xi\rangle
+2\langle JX,Z\rangle\langle Y,J\xi\rangle\right).
\]

\item\label{thCodazzi3Spaces} Let $X\in\Gamma(T_\alpha)$,
    $Y\in\Gamma(T_\beta)$ and $Z\in\Gamma(T_\gamma)$. Then
\[
\langle\bar{R}(X,Y)Z,\xi\rangle
=(\beta-\gamma)\langle\nabla_X Y,Z\rangle
-(\alpha-\gamma)\langle\nabla_Y X,Z\rangle.
\]
\end{enumerate}
\end{lemma}

The Gauss equation implies (again, see \cite[Lemma 4]{BD05} for a proof)

\begin{lemma}\label{thGauss}
Let $X\in\Gamma(T_\alpha)$ and $Y\in\Gamma(T_\beta)$, with
$\alpha\neq\beta$, be unit vector fields. Then
\begin{eqnarray*}
0 &=& (\beta-\alpha)(-c-4\alpha\beta-2c\langle JX,Y\rangle^2
+8\langle\nabla_X Y,\nabla_Y X\rangle-4\langle\nabla_X X,\nabla_Y Y\rangle)\\
&&-4c\langle JX,Y\rangle(X\langle Y,J\xi\rangle+Y\langle X,J\xi\rangle)\\
&&-c\langle X,J\xi\rangle(3Y\langle JX,Y\rangle+\langle\nabla_Y X,JY\rangle
-2\langle\nabla_X Y,JY\rangle)\\
&&-c\langle Y,J\xi\rangle(3X\langle JX,Y\rangle-\langle\nabla_X Y,JX\rangle
+2\langle\nabla_Y X,JX\rangle).
\end{eqnarray*}
\end{lemma}

%%%%%%%%%%%%%%%%%%%%%%%%%% Subsection %%%%%%%%%%%%%%%%%%%%%%%%%%%%%%%
\subsection{Discussion of examples}\label{secExamples}

Part (\ref{thMainCPn}) of the Main Theorem states that there are no examples
of real hypersurfaces with constant principal curvatures in $\C P^n(c)$ whose
Hopf vector field has $h=2$ nontrivial projections onto the principal
curvature spaces of $M$. Thus, we will focus on describing briefly the
examples of part (\ref{thMainCHn}) of the Main Theorem. These examples where
first constructed in \cite{BB01} and their geometry was studied in
\cite{BD09}. %See also \cite{HHT09}.

%The connected Lie group $G=SU(1,n)$ acts transitively on $\C H^n(c)$. Fix
%a point $o\in\C H^n(c)$ and let $K$ be the isotropy group of $G$ at $o$.
%Then $K$ is isomorphic to $S(U(1)U(n))$ and $(G,K)$ is a symmetric pair.
%Then, $\C H^n(c)$ may be identified with the quotient $G/K$. The group
%$G$ is simple and of noncompact type. Write $\g{g}$ for the Lie algebra
%of $G$ and $\g{k}$ for the Lie algebra of $K$. Let $\g{g}=\g{k}\oplus\g{p}$
%be the Cartan decomposition of $\g{g}$ with respect to $o\in\C H^n(c)$ and
%let $\theta$ be the corresponding Cartan involution. If $B$ denotes the
%Killing form of $\g{g}$ then we define the inner product $B_\theta$ on $\g{g}$
%by $B_\theta(X,Y)=-B(X,\theta Y)$ for all $X,Y\in\g{g}$. If we identify
%$\g{p}\cong T_o\C H^n(c)$ as usual, then the inner product of
%$\C H^n(c)\cong G/K$ is the left-invariant inner product determined by
%$\langle X,Y\rangle=-\frac{1}{c(n+1)}B_\theta(X,Y)$ for all $X,Y\in\g{p}$
%(this choice gives precisely the complex hyperbolic space of constant
%holomorphic sectional curvature $c<0$).

The connected simple Lie group $G=SU(1,n)$ acts transitively on $\C
H^n(c)$. Fix a point $o\in\C H^n(c)$ and let $K$ be the
isotropy group of $G$ at $o$. The subgroup $K$ of $G$ is isomorphic to
$S(U(1)U(n))$. Furthermore, $(G,K)$ is a symmetric pair and $\C
H^n(c)$ may be identified with the quotient $G/K$. Write $\g{g}$ for the Lie
algebra of $G$ and $\g{k}$ for the Lie algebra of $K$. Let
$\g{g}=\g{k}\oplus\g{p}$ be the Cartan decomposition of $\g{g}$
with respect to $o\in\C H^n(c)$. We choose a
maximal abelian subspace $\g{a}$ of $\g{p}$; then, $\dim\g{a}=1$
since $\C H^n(c)$ has rank one. Let
$\g{g}=\g{g}_{-2\alpha}\oplus\g{g}_{-\alpha}\oplus\g{g}_0
\oplus\g{g}_\alpha\oplus\g{g}_{2\alpha}$ be the restricted root space
decomposition of $\g{g}$ with respect to $\g{a}$ and assume
that $\alpha$ is a positive root. Then,
$\g{n}=\g{g}_\alpha\oplus\g{g}_{2\alpha}$ is a 2-step nilpotent
subalgebra of $\g{g}$ isomorphic to the
$(2n-1)$-dimensional Heisenberg algebra. Furthermore,
$\g{g}=\g{k}\oplus\g{a}\oplus\g{n}$ is an Iwasawa decomposition
of $\g{g}$. If $A$ and $N$ denote the connected subgroups
of $G$ whose Lie algebras are $\g{a}$ and $\g{n}$,
then $G=KAN$ is an Iwasawa decomposition of $G$. The solvable
group $AN$ is simply connected and acts simply transitively on
$\C H^n(c)$. Thus, we can identify $\g{a}\oplus\g{n}$ with
$T_o\C H^n(c)$. The Riemannian metric of $\C H^n(c)$ induces a
left-invariant metric on $AN$ which makes $AN$ isometric to $\C H^n(c)$. Similarly, the complex structure
$J$ on $T_o\C H^n(c)$ induces a complex structure on
$\g{a}\oplus\g{n}$ which we also denote by $J$. We have $J\g{a}=\g{g}_{2\alpha}$, and $\g{g}_\alpha$ is
$J$-invariant. Let $B\in\g{a}$ be a unit vector and define
$Z=JB\in\g{g}_{2\alpha}$.

Let $\g{w}$ be a linear subspace of $\g{g}_\alpha$ such that
the orthogonal complement
$\g{w}^\perp=\g{g}_\alpha\ominus\g{w}$ of $\g{w}$ in
$\g{g}_\alpha$ has constant K\"{a}hler angle $\varphi\in(0,\pi/2]$,
that is, the angle between $Jv$ and $\g{w}^\perp$ is $\varphi$
for all nonzero $v\in\g{w}^\perp$. Then, $\varphi=\pi/2$ if and only if
$\g{w}^\perp$ is real, or equivalently, if and only if $J\g{w}^\perp$ is orthogonal to $\g{w}^\perp$. Let $k$ be the dimension of
$\g{w}^\perp$. Then,
$\g{s}=\g{a}\oplus\g{w}\oplus\g{g}_{2\alpha}$ is a subalgebra
of $\g{a}\oplus\g{n}$. Let $S$ be the connected simply
connected subgroup of $AN$ whose Lie algebra is $\g{s}$. We
define the Berndt-Br\"{u}ck submanifolds as \cite{BB01} (see \cite{L99} for $k=1$)
\[
W^{2n-k}_\varphi=S\cdot o,\quad\mbox{ and }\quad W^{2n-k}=W^{2n-k}_{\pi/2}.
\]
The Berndt-Br\"{u}ck submanifolds $W^{2n-k}_\varphi$ are
homogeneous, have normal bundle of rank $k$ and constant K\"{a}hler
angle $\varphi\in(0,\pi/2]$, and their second fundamental form
$\II$ is given by the trivial symmetric bilinear extension of
$\II(Z,P\xi)=(\sin(\varphi)\,\sqrt{-c}/2)\xi$ for all
$\xi\in\g{w}^\perp$, where $P\xi$ is the orthogonal projection
of $J\xi$ onto $TW_\varphi^{2n-k}$. In particular, the
submanifolds $W^{2n-k}_\varphi$ are minimal, and ruled by the
totally geodesic complex hyperbolic subspaces determined by
their maximal holomorphic distribution. If $\varphi=\pi/2$ then
$P=J$ and the Berndt-Br\"{u}ck submanifolds have totally real
normal bundle. Conversely \cite[Theorem 1]{BD09}

\begin{theorem}\label{thRigidity}
Let $M$ be a $(2n-k)$-dimensional connected submanifold in $\C
H^n(c)$, $n\geq 2$, with normal bundle $\nu
M$ of constant K\"{a}hler angle $\varphi\in(0,\pi/2]$. Assume that there exists a unit vector
field $Z$ tangent to the maximal holomorphic subbundle of
$TM$ such that the second fundamental form
$\II$ of $M$ is given by the trivial symmetric bilinear extension of
\[
\II(Z,P\xi)=\sin(\varphi)\,\frac{\sqrt{-c}}{2}\,\xi
\]
for all $\xi\in\Gamma(\nu M)$. Then $M$ is holomorphically
congruent to an open part of the ruled minimal submanifold
$W^{2n-k}_\varphi$.
\end{theorem}

In particular, the Berndt-Br\"{u}ck submanifolds $W^{2n-k}$ are
determined by the equation $\II(Z,J\xi)=({\sqrt{-c}}/{2})\xi$
and the fact that their normal bundle is totally real. Geometrically, they are constructed in the following way. Fix a horosphere $\mathcal{H}$ in a totally geodesic real hyperbolic space $\R H^{k+1}(c)\subset\C H^n(c)$. Attach at each point the totally geodesic $\C H^{n-k}(c)$ which is tangent to the orthogonal complement of the complex span of the tangent space of $\mathcal{H}$ at $p$. The resulting submanifold is congruent to $W^{2n-k}$.

Let $N_K^0(S)$ denote the connected component of the identity
transformation of the normalizer of $S$ in $K$. Then,
$N_K^0(S)S$ acts on $\C H^n(c)$ with cohomogeneity one and
$W^{2n-k}_\varphi=N_K^0(S)S\cdot o$. If $k>1$, then the
principal orbits of $N_K^0(S)S$ are tubes around
$W^{2n-k}_\varphi$. If $k=1$, then $\varphi=\pi/2$, the
action of $N_K^0(S)S$ is orbit equivalent to the action of $S$,
and its orbits form a homogeneous foliation on $\C H^n(c)$ that
was first studied in \cite{Be98}.

Let $M$ be a principal orbit of $N_K^0(S)S$. If
$\varphi\in(0,\pi/2)$ then the Hopf vector field of $M$ has
$h=3$ nontrivial projections onto the principal curvature
spaces of $M$. If $\varphi=\pi/2$, then the Hopf vector field
of $M$ has $h=2$ nontrivial projections onto the principal
curvature spaces of $M$. The objective of
part~(\ref{thMainCHn}) of the Main Theorem is to give a
geometric characterization of the tubes around $W^{2n-k}$,
$k\in\{2,\dots,n-1\}$, and the equidistant hypersurfaces to
$W^{2n-1}$.

%%%%%%%%%%%%%%%%%%%%%%%%%% Section %%%%%%%%%%%%%%%%%%%%%%%%%%%%%%%
\section{Proof of the Main Theorem}\label{secProof}

In this section we prove the Main Theorem. Our main goal is to describe accurately the eigenvalue structure of a real hypersurface in the conditions of the Main Theorem (Theorem~\ref{thSummary}). Then we finish the
proof using standard Jacobi field theory (\S\ref{secJacobi}).
%We start we some algebraic preliminaries.
%Some
%steps of our proof will be similar to some steps in the proof
%given in \cite{BD05}. When this is the case we will point it
%out and refer to \cite{BD05} for more details. However, several
%parts of the proof are fundamentally different, so we will
%explain the new techniques in a greater detail.

%%%%%%%%%%%%%%%%%%%%%%%%%% Subsection %%%%%%%%%%%%%%%%%%%%%%%%%%%%%%%
\subsection{Notation and setup}\label{secNotation}

Let $M$ be a connected real hypersurface with $g>1$ distinct constant
principal curvatures in a complex space form $\bar{M}(c)$.
Since the calculations that follow are local we may assume that
we have a globally defined unit normal vector field $\xi$. We
denote by $\lambda_1,\dots,\lambda_g$ the principal curvatures of $M$.

By assumption, the number of nontrivial projections of
$J\xi$ onto the principal curvature distributions
$T_{\lambda_i}$, $i\in\{1,\dots,g\}$, is $h=2$. By relabeling the indices we may also assume
that $J\xi$ has nontrivial projection onto
$T_{\lambda_1}$ and $T_{\lambda_2}$. Hence, there exist unit
vectors fields $U_i\in\Gamma(T_{\lambda_i})$, $i\in\{1,2\}$, and
positive smooth functions $b_i:M\to\R$,
$i\in\{1,2\}$, such that
\[
J\xi=b_1 U_1+b_2 U_2.
\]
Obviously, $b_1^2+b_2^2=1$. Moreover,

\begin{lemma}\label{thA}
We have $g\geq 3$, $\langle JU_1,U_2\rangle=0$ and there exists
a unit vector field $A\in\Gamma(\oplus_{k=3}^g T_{\lambda_k})$
such that
\begin{eqnarray*}
JU_i    &=& (-1)^i b_j A-b_i\xi, \quad(i,j\in\{1,2\}, i\neq j),\\
JA      &=& b_2 U_1-b_1 U_2.
\end{eqnarray*}
\end{lemma}

\begin{proof}
The proof is similar to that of \cite[Lemma 7]{BD05}, so we
just sketch it. We will assume in what follows $i,j\in\{1,2\}$,
$i\neq j$, and $k\in\{3,\dots,g\}$.

Since $T_{\lambda_i}$, $i\in\{1,2\}$, is real by
Lemma~\ref{thCodazzi}~(\ref{thCodazzi1Space}) we can write
$JU_i=\langle JU_i,U_j\rangle U_j+W_{ij}+\sum_{k=3}^g
W_{ik}-b_i\xi$, where $W_{ij}\in\Gamma(T_{\lambda_j}\ominus\R
U_j)$ %(orthogonal complement of $\R U_j$ in $T_{\lambda_j}$),
and $W_{ik}\in\Gamma(T_{\lambda_k})$. (Here and henceforth, the symbol $\ominus$ is used to denote orthogonal complement.) From $J\xi=b_1 U_1+b_2 U_2$
we get
\[
-\xi=J^2\xi=b_2(\langle JU_2,U_1\rangle U_1+W_{21})
+b_1(\langle JU_1,U_2\rangle U_2+W_{12})
+\sum_{k=3}^g(b_1W_{1k}+b_2W_{2k})-\xi.
\]
Thus, $g\geq 3$, $\langle JU_1,U_2\rangle=0$, $W_{12}=W_{21}=0$,
and $b_1W_{1k}+b_2W_{2k}=0$ for all $k$. If we define $A\in\Gamma(\oplus_{k=3}^g T_{\lambda_k})$ by $\sum_{k=3}^gW_{ik}=(-1)^i b_j A$, then the last equality implies $\sum_{k=3}^gW_{jk}=(-1)^j b_i A$ (recall $i,j\in\{1,2\}$,
$i\neq j$). This gives the desired expression for $JU_i$, $i\in\{1,2\}$. Finally, from $b_1^2+b_2^2=1$ and
$-U_1=J(JU_1)=-b_2JA-b_1J\xi=-b_2JA-U_1+b_2^2U_1-b_1b_2U_2$
we obtain $JA=b_2U_1-b_1U_2$.
\end{proof}

%In particular, the distribution $\R\xi\oplus\R U_1\oplus\R U_2\oplus\R A$ of
%$T\bar{M}(c)$ is complex.

%%%%%%%%%%%%%%%%%%%%%%%%%% Subsection %%%%%%%%%%%%%%%%%%%%%%%%%%%%%%%
\subsection{The vector field $A$}\label{secA}

In view of Lemma~\ref{thA} we may write
\[
A=\sum_{k=3}^g A_k,\mbox{ with } A_k\in\Gamma(T_{\lambda_k}), k\in\{3,\dots g\}.
\]
The aim of this subsection is to show that all but one $A_k$ are
zero and hence we can assume for example that
$A\in\Gamma(T_{\lambda_3})$ (Proposition~\ref{thAinT3}).
%Note that in \cite{BD05} the
%assumption $g=3$ already implied this property. However, here
%$g$ is undetermined and thus we need an additional argument.
The main difficulty here is the fact that $g$ is not known.
We start with the following

\begin{lemma}\label{thNablaUiUj}
Let $i,j\in\{1,2\}$ with $i\neq j$. Then we have
%\begin{eqnarray*}
%\nabla_{U_i}U_i&=&
%    \sum_{k=3}^g (-1)^j\frac{3cb_1b_2}{4(\lambda_k-\lambda_i)}A_k,\\
%\nabla_{U_i}U_j&=&
%    \sum_{k=3}^g (-1)^j\left(\lambda_i
%        -\frac{3cb_i^2}{4(\lambda_k-\lambda_i)}\right)A_k.\\
%\end{eqnarray*}
\[
\nabla_{U_i}U_i=
    \sum_{k=3}^g (-1)^j\frac{3cb_1b_2}{4(\lambda_k-\lambda_i)}A_k,\qquad
\nabla_{U_i}U_j=
    \sum_{k=3}^g (-1)^j\left(\lambda_i
        -\frac{3cb_i^2}{4(\lambda_k-\lambda_i)}\right)A_k.
\]
\end{lemma}

\begin{proof}
Again, this is quite similar to \cite[Lemma 8]{BD05}. We assume $i,j\in\{1,2\}$, $i\neq j$, and
$k\in\{3,\dots,g\}$. Let $W_i\in\Gamma(T_{\lambda_i}\ominus\R
U_i)$ and $W_k\in\Gamma(T_{\lambda_k}\ominus\R A_k)$.

Since $U_i$ has unit length we get
$\langle\nabla_{U_i}U_i,U_i\rangle=0$.
Lemma~\ref{thCodazzi}~(\ref{thCodazzi2Spaces}) yields
$\langle\nabla_{U_i}U_i,U_j\rangle
=\langle\nabla_{U_i}U_i,W_j\rangle=\langle\nabla_{U_i}U_i,W_k\rangle=0$,
and $\langle\nabla_{U_i}U_i,A_k\rangle =3(-1)^j
cb_1b_2/(4(\lambda_k-\lambda_i))$.
%\[
%\langle\nabla_{U_i}U_i,U_j\rangle
%=\langle\nabla_{U_i}U_i,W_j\rangle=\langle\nabla_{U_i}U_i,W_k\rangle=0,\quad
%\langle\nabla_{U_i}U_i,A_k\rangle
%=(-1)^j\frac{3cb_1b_2}{4(\lambda_k-\lambda_i)}.
%\]
From $\bar\nabla J=0$, the Weingarten formula, and
Lemma~\ref{thA}, we obtain $\langle W_i,\bar\nabla_{U_i}J\xi\rangle=-\lambda_i\langle W_i,JU_i\rangle=0$. Hence, using $J\xi=b_1U_1+b_2U_2$, and Lemma~\ref{thCodazzi}~(\ref{thCodazzi2Spaces}),  we get
\[
0=U_i\langle W_i,J\xi\rangle
=\langle\nabla_{U_i}W_i,J\xi\rangle+\langle W_i,\bar\nabla_{U_i}J\xi\rangle
=-b_i\langle\nabla_{U_i}U_i,W_i\rangle.
\]
Since $b_i\neq 0$ the expression for $\nabla_{U_i}U_i$ follows.

As $U_j$ has unit length,
$\langle\nabla_{U_i}U_j,U_j\rangle=0$. From
Lemma~\ref{thCodazzi}~(\ref{thCodazzi2Spaces}) we obtain
$\langle\nabla_{U_i}U_j,U_i\rangle=\langle\nabla_{U_i}U_j,W_i\rangle=0$. Now, the Weingarten formula and Lemma~\ref{thA} imply $\langle W_j,\bar\nabla_{U_i}J\xi\rangle=-\lambda_i\langle W_j,J U_i\rangle=0$, and thus, Lemma~\ref{thCodazzi}~(\ref{thCodazzi2Spaces}), yields
\[
0=U_i\langle W_j,J\xi\rangle
=\langle\nabla_{U_i}W_j,J\xi\rangle+\langle W_j,\bar\nabla_{U_i}J\xi\rangle
=b_j\langle\nabla_{U_i}W_j,U_j\rangle.
\]
This implies $\langle\nabla_{U_i}W_j,U_j\rangle=0$. A similar calculation gives $\langle\nabla_{U_i}W_k,U_j\rangle=0$. Finally, by Lemma~\ref{thCodazzi}~(\ref{thCodazzi2Spaces}), and Lemma~\ref{thA} we have
\[
0=U_i \langle A_k,J\xi\rangle
=\langle\nabla_{U_i}A_k,J\xi\rangle+\langle A_k,\bar\nabla_{U_i}J\xi\rangle
=(-1)^i\frac{3c b_i^2 b_j}{4(\lambda_k-\lambda_i)}-b_j\langle\nabla_{U_i}U_j,A_k\rangle
-(-1)^i\lambda_i b_j,
\]
from where we get $\langle\nabla_{U_i}U_j,A_k\rangle$.
Altogether this yields the formula for $\nabla_{U_i}U_j$.
\end{proof}

%\begin{corollary}%
%Let $k\in\{3,\dots,g\}$ and $p\in M$ such that $(A_k)_p\neq 0$. Then
%\[
%\frac{3c(\lambda_2-\lambda_k)}{4(\lambda_1-\lambda_k)}b_1(p)^2
%+\frac{3c(\lambda_1-\lambda_k)}{4(\lambda_2-\lambda_k)}b_2(p)^2
%=-\frac{c}{4}-\lambda_1(\lambda_2-\lambda_k)-\lambda_2(\lambda_1-\lambda_k).
%\]
%\end{corollary}

Now we can prove the main result of this section.

\begin{proposition}\label{thAinT3}
$A\in\Gamma(T_{\lambda_k})$ for some $k\in\{3,\dots,g\}$.
\end{proposition}

\begin{proof}
On the contrary, assume that there exists a point $p\in M$ and
two distinct integers $r,s\in\{3,\dots,g\}$ such that
$(A_r)_p,(A_s)_p\neq 0$. Hence, in a neighborhood of $p$ we have
$A_r,A_s\neq 0 $ as well. We will work in that neighborhood
from now on.

Applying Lemma~\ref{thCodazzi}~(\ref{thCodazzi3Spaces}) to the vector fields
$U_1$, $U_2$, and $A_k$, $k\in\{r,s\}$, and using Lemma~\ref{thNablaUiUj} we easily get
\begin{equation}\label{eqCodazziU1U2Ak}
\frac{3c(\lambda_2-\lambda_k)}{4(\lambda_1-\lambda_k)}b_1^2
+\frac{3c(\lambda_1-\lambda_k)}{4(\lambda_2-\lambda_k)}b_2^2
=-\frac{c}{4}-\lambda_1(\lambda_2-\lambda_k)-\lambda_2(\lambda_1-\lambda_k),
\quad k\in\{r,s\}.
\end{equation}
Together with $b_1^2+b_2^2=1$, this
yields a linear system of three equations with unknowns
$b_1^2$ and $b_2^2$. This system must be compatible. We show it is
determined (that is, it has a unique solution). If it were not,
the rank of the system would, at most, be one. In particular,
\[
\left\lvert\begin{array}{cc}
\frac{3c(\lambda_2-\lambda_k)}{4(\lambda_1-\lambda_k)}
& \frac{3c(\lambda_1-\lambda_k)}{4(\lambda_2-\lambda_k)}\\
1 & 1
\end{array}\right\rvert=
3c\frac{(\lambda_2-\lambda_1)(\lambda_1+\lambda_2-2\lambda_k)}
{4(\lambda_1-\lambda_k)(\lambda_2-\lambda_k)}=0,\quad k\in\{r,s\},
\]
which implies $\lambda_1+\lambda_2-2\lambda_k=0$, $k\in\{r,s\}$, and
hence $\lambda_r=\lambda_s$, contradiction. We conclude that
the above system is determined. Therefore, we can find an
expression for $b_1^2$ and $b_2^2$ in terms of the principal
curvatures and $c$. Since these are constant, it follows that
$b_1$ and $b_2$ are constant.

We take $i,j\in\{1,2\}$, $i\neq j$, and $k\in\{r,s\}$. Since
$b_i$ is constant and $U_i$ has unit length, using $J\xi=b_1U_1+b_2U_2$, the
Weingarten formula, and Lemma~\ref{thA} we get
\[
0=A_k(b_i)=A_k\langle U_i,J\xi\rangle
=\langle\nabla_{A_k}U_i,J\xi\rangle+\langle U_i,\bar\nabla_{A_k}J\xi\rangle
=b_j\langle\nabla_{A_k}U_i,U_j\rangle-(-1)^j b_j\lambda_k,
\]
and thus, $\langle\nabla_{A_k}U_i,U_j\rangle=(-1)^j\lambda_k$.
Taking this, Lemma~\ref{thA}, and Lemma~\ref{thNablaUiUj} into account,
Lemma~\ref{thCodazzi}~(\ref{thCodazzi3Spaces}) for $A_k$, $U_1$
and $U_2$ yields
\[
\frac{c}{4}(2b_2^2-b_1^2)=\langle\bar{R}(A_k,U_1)U_2,\xi\rangle
=(\lambda_1-\lambda_2)\lambda_k+(\lambda_k-\lambda_2)\left(\lambda_1
-\frac{3cb_1^2}{4(\lambda_k-\lambda_1)}\right),\quad k\in\{r,s\}.
\]
We can rearrange this as:
\begin{equation}\label{eqCodazziAkU1U2}
\left(\frac{c}{4}-\frac{3c(\lambda_k-\lambda_2)}{4(\lambda_k-\lambda_1)}\right)b_1^2
-\frac{c}{2}b_2^2
=(\lambda_2-\lambda_1)\lambda_k+\lambda_1(\lambda_2-\lambda_k),
\quad k\in\{r,s\}.
\end{equation}
Hence, \eqref{eqCodazziU1U2Ak}, \eqref{eqCodazziAkU1U2}, and $b_1^2+b_2^2=1$ give a linear system of five equations with unknowns $b_1^2$ and
$b_2^2$.
This system is compatible by assumption, so it has rank two.
Then, all minors of order three of the augmented matrix of the system vanish. This implies (take \eqref{eqCodazziU1U2Ak}, \eqref{eqCodazziAkU1U2}, and $b_1^2+b_2^2=1$, with $k\in\{r,s\}$, and then both equations in \eqref{eqCodazziAkU1U2} and $b_1^2+b_2^2=1$):
\begin{eqnarray}
\frac{3c(\lambda_1-\lambda_2)^2(-12\lambda_k^2+8\lambda_1\lambda_k
+8\lambda_2\lambda_k+c-4\lambda_1\lambda_2)}
{16(\lambda_1-\lambda_k)(\lambda_k-\lambda_2)}
&=&0,\ k\in\{r,s\},\label{eqlambdars1}\\
\frac{3c(\lambda_2-\lambda_1)(\lambda_r-\lambda_s)
(4\lambda_1^2-4\lambda_r\lambda_1-4\lambda_s\lambda_1
+c+2\lambda_2\lambda_r+2\lambda_2\lambda_s)}
{8(\lambda_1-\lambda_r)(\lambda_1-\lambda_s)}   &=& 0.\label{eqlambdars2}
\end{eqnarray}
In particular, \eqref{eqlambdars1} implies $-12\lambda_k^2+8\lambda_1\lambda_k
+8\lambda_2\lambda_k+c-4\lambda_1\lambda_2=0$. Putting $k=r$ and $k=s$, and subtracting, we get
$4(2\lambda_1+2\lambda_2-3\lambda_r-3\lambda_s)(\lambda_r-\lambda_s)=0$,
from where we obtain
$\lambda_r+\lambda_s=2(\lambda_1+\lambda_2)/3$. Taking this
into account, \eqref{eqlambdars2} gives
$(4\lambda_1^2-4\lambda_1\lambda_2+4\lambda_2^2+3c)/3=0$. The
discriminant of $-12\lambda_k^2+8\lambda_1\lambda_k
+8\lambda_2\lambda_k+c-4\lambda_1\lambda_2=0$ as a quadratic equation in
$\lambda_k$ is precisely
$16(4\lambda_1^2-4\lambda_1\lambda_2+4\lambda_2^2+3c)$, so this
discriminant vanishes. As a consequence, this quadratic equation
has a unique solution and hence $\lambda_r=\lambda_s$. This is
a contradiction. Therefore, all but one $A_k$,
$k\in\{3,\dots,g\}$, are zero for each $p$. The result follows
by continuity.
\end{proof}

%%%%%%%%%%%%%%%%%%%%%%%%%% Subsection %%%%%%%%%%%%%%%%%%%%%%%%%%%%%%%
\subsection{Some properties of the principal curvature
spaces}\label{secPrincipalSpaces}

In view of Proposition \ref{thAinT3}, we may assume from now on
that $A\in\Gamma(T_{\lambda_3})$. Moreover, we can choose an
orientation on $M$ and a relabeling of the indices so that
\[
\lambda_1<\lambda_2,\quad\mbox{ and }\quad\lambda_3\geq 0.
\]
We will follow this convention from now on.

First we calculate some covariant derivatives.

\begin{lemma}\label{thNablaUiUjA}
Let $i,j\in\{1,2\}$ with $i\neq j$. Then we have
\begin{eqnarray}
\nabla_{U_i}U_i &=& (-1)^j\frac{3cb_1b_2}
    {4(\lambda_3-\lambda_i)}A,\label{eqUiUi}\\
\nabla_{U_i}U_j &=& (-1)^j\left(\lambda_i-\frac{3cb_i^2}
    {4(\lambda_3-\lambda_i)}\right)A,\label{eqUiUj}\\
\nabla_{U_i}A   &=& (-1)^i\frac{3cb_1b_2}{4(\lambda_3-\lambda_i)}U_i
    +(-1)^i\left(\lambda_i-\frac{3cb_i^2}
    {4(\lambda_3-\lambda_i)}\right)U_j,\label{eqUiA}\\
\nabla_{A}U_i &=& \frac{(-1)^j}{\lambda_i-\lambda_j}
    \left(\frac{c(2b_j^2-b_i^2)}{4}+(\lambda_j-\lambda_3)
    \left(\lambda_i-\frac{3cb_i^2}
    {4(\lambda_3-\lambda_i)}\right)\right)U_j,\label{eqAUi}\\[2ex]
\nabla_{A}A &=& 0.\label{eqAA}
\end{eqnarray}
\end{lemma}

\begin{proof}
The proof is similar to that of \cite[Lemma 8]{BD05}. Equations \eqref{eqUiUi} and \eqref{eqUiUj} are a
direct consequence of Lemma~\ref{thNablaUiUj} and Proposition~\ref{thAinT3}. Assume $i,j\in\{1,2\}$, $i\neq j$, and $k\in\{4,\dots,g\}$. Let
$W_i\in\Gamma(T_{\lambda_i}\ominus\R U_i)$,
$W_3\in\Gamma(T_{\lambda_3}\ominus\R A)$ and
$W_k\in\Gamma(T_{\lambda_k})$.

According to \eqref{eqUiUi} and \eqref{eqUiUj}, in order to
prove \eqref{eqUiA} we have to show
$\langle\nabla_{U_i}A,A\rangle=0$ (obvious because $A$ is a
unit vector field), and $\langle\nabla_{U_i}A,W_l\rangle=0$ for all
$l\in\{1,\dots,g\}$. The latter follows after using
$\bar\nabla J=0$, the
Weingarten formula, Lemma~\ref{thA}, and \eqref{eqUiUi}, with
\begin{eqnarray*}
0&=&
U_i\langle JU_i,W_l\rangle\
=\ \langle\bar\nabla_{U_i}JU_i,W_l\rangle
+\langle JU_i,\bar\nabla_{U_i}W_l\rangle\\
&=&-\langle\nabla_{U_i}U_i,JW_l\rangle
%+\lambda_i\langle \xi,JW_l\rangle
+(-1)^ib_j\langle A,\nabla_{U_i}W_l\rangle
-b_i\langle \xi,\bar\nabla_{U_i}W_l\rangle
=(-1)^j b_j\langle\nabla_{U_i}A,W_l\rangle.
\end{eqnarray*}

We now prove \eqref{eqAUi}. Obviously, $\langle\nabla_A
U_i,U_i\rangle=0$, and $\langle\nabla_A U_i,A\rangle=0$ by
Lemma~\ref{thCodazzi}~(\ref{thCodazzi2Spaces}). Applying
Lemma~\ref{thCodazzi}~(\ref{thCodazzi3Spaces}) to $A$, $U_i$
and $U_j$, using Lemma~\ref{thA} and \eqref{eqUiUj}, gives
\[
\frac{c}{4}(-1)^i(b_i^2-2b_j^2)
=(\lambda_i-\lambda_j)\langle\nabla_A U_i,U_j\rangle
-(\lambda_3-\lambda_j)(-1)^i\left(\lambda_i
-\frac{3cb_i^2}{4(\lambda_3-\lambda_i)}\right),
\]
from where we get $\langle\nabla_A U_i,U_j\rangle$. For
$l\in\{j,3,\dots,g\}$, a similar argument with
Lemma~\ref{thCodazzi}~(\ref{thCodazzi3Spaces}) applied to $A$, $U_i$,
and $W_l$, taking Lemma~\ref{thA} and \eqref{eqUiA} into account, yields
$\langle\nabla_A{U_i},W_l\rangle=0$. Finally, the previous equality (interchanging $i$ and $j$ and putting $l=i$) gives
\begin{eqnarray*}
0&=&
A\langle W_i,J\xi\rangle\
=\ \langle\nabla_A W_i,J\xi\rangle+\langle W_i,\bar\nabla_A J\xi\rangle\\
&=&b_i\langle\nabla_A W_i,U_i\rangle+b_j\langle\nabla_A W_i,U_j\rangle
-\lambda_3\langle W_i,JA\rangle
=-b_i\langle\nabla_A U_i,W_i\rangle.
\end{eqnarray*}
Altogether this proves \eqref{eqAUi}.

We have $\langle\nabla_A A,A\rangle=0$, and $\langle\nabla_A
A,U_i\rangle=\langle\nabla_A A,W_l\rangle=0$ for all
$l\in\{1,2,4,\dots,g\}$ by
Lemma~\ref{thCodazzi}~(\ref{thCodazzi2Spaces}). From
$\bar\nabla J=0$, \eqref{eqAUi}, Lemma~\ref{thA}, and the
Weingarten formula we get
\begin{eqnarray*}
0&=&
A\langle JU_i,W_3\rangle\
=\ \langle\bar\nabla_{A}JU_i,W_3\rangle
+\langle JU_i,\bar\nabla_{A}W_3\rangle\\
&=&-\langle\nabla_{A}U_i,JW_3\rangle
+(-1)^ib_j\langle A,\nabla_{A}W_3\rangle
-b_i\langle\xi,\bar\nabla_{A}W_3\rangle
=(-1)^j b_j\langle\nabla_{A}A,W_3\rangle.
\end{eqnarray*}
from where \eqref{eqAA} follows.
\end{proof}

Our main difficulty from now on is the fact that the number
$g$ of principal curvatures is not known. In fact, the aim of
Subsection~\ref{secBoundg} is to obtain a bound on $g$.
%From now on and until Theorem~\ref{thSummary} the method we use
%to prove our Main Theorem is different from that of
%\cite{BD05}.
An important step in the proof is the following

\begin{proposition}\label{thB}
The functions $b_1$ and $b_2$ are constant. In fact
\[
b_i^2=\frac{4(\lambda_j-2\lambda_3)(\lambda_i-\lambda_3)^2}
{c(\lambda_i-\lambda_j)},\quad (i,j\in\{1,2\},i\neq j).
\]
Moreover,
$c-4\lambda_1\lambda_2+8(\lambda_1+\lambda_2)\lambda_3-12\lambda_3^2=0$.
\end{proposition}

\begin{proof}
First we show that the functions $b_1$ and $b_2$ are constant.

We apply Lemma~\ref{thGauss} to $U_1$ and $U_2$, using
Lemma~\ref{thA} and Lemma~\ref{thNablaUiUjA},
\begin{eqnarray*}
0\!\!&=&
(\lambda_2-\lambda_1)(-c-4\lambda_1\lambda_2
+8\langle\nabla_{U_1}U_2,\nabla_{U_2}U_1\rangle
-4\langle\nabla_{U_1}U_1,\nabla_{U_2}U_2\rangle)\\
&&-cb_1(\langle\nabla_{U_2}U_1,JU_2\rangle-2\langle\nabla_{U_1}U_2,JU_2\rangle)
-cb_2(-\langle\nabla_{U_1}U_2,JU_1\rangle+2\langle\nabla_{U_2}U_1,JU_1\rangle)\\
&=&
-(\lambda_2-\lambda_1)(c+12\lambda_1\lambda_2)-\frac{3c^2}{2(\lambda_3-\lambda_1)}b_1^4
+\frac{3c^2}{2(\lambda_3-\lambda_2)}b_2^4
+\frac{3c^2(\lambda_1-\lambda_2)}{(\lambda_3-\lambda_1)(\lambda_3-\lambda_2)}
b_1^2b_2^2\\
&&+\frac{c(6\lambda_2^2-7\lambda_1\lambda_2-2\lambda_1^2+2\lambda_1\lambda_3
+\lambda_2\lambda_3)}{\lambda_3-\lambda_1}b_1^2
-\frac{c(6\lambda_1^2-7\lambda_1\lambda_2-2\lambda_2^2+2\lambda_2\lambda_3
+\lambda_1\lambda_3)}{\lambda_3-\lambda_2}b_2^2.
\end{eqnarray*}
Now we substitute $b_2^2$ by $1-b_1^2$ to get
\[
0=\frac{9c^2(\lambda_2-\lambda_1)}
{2(\lambda_3-\lambda_1)(\lambda_3-\lambda_2)}b_1^4+\Lambda_1b_1^2+\Lambda_0,
\]
where $\Lambda_1$ and $\Lambda_0$ are constants depending on
$c$, $\lambda_1$, $\lambda_2$, and $\lambda_3$. This equation is
a quadratic equation in $b_1^2$ and the coefficient of $b_1^4$
does not vanish. Hence, it has at most two real solutions
depending on the constants $c$, $\lambda_1$, $\lambda_2$ and
$\lambda_3$. Since $M$ is connected it follows that $b_1$ and
$b_2$ are constant.

From the argument above one might derive an explicit
expression for $b_i$, $i\in\{1,2\}$. However, that expression
would involve square roots that would make later calculations
difficult. Instead, we use the constancy of these functions to
give an alternative formula which is easier to handle. For
$i\in\{1,2\}$, using lemmas~\ref{thA} and~\ref{thNablaUiUjA}, and
the Weingarten formula, we get
\begin{eqnarray*}
0&=&
A(b_i)\ =\ A\langle U_i,J\xi\rangle\ =\
\langle\nabla_A U_i,J\xi\rangle+\langle U_i,\bar\nabla_A J\xi\rangle\
=\ b_j\langle\nabla_A U_i,U_j\rangle-\lambda_3\langle U_i,JA\rangle\\
&=&
(-1)^ib_j\left(c\frac{-\lambda_i+3\lambda_j-2\lambda_3}
{4(\lambda_i-\lambda_j)(\lambda_3-\lambda_i)}b_i^2
-\frac{c}{2(\lambda_i-\lambda_j)}b_j^2
+\frac{2\lambda_i\lambda_3-\lambda_j\lambda_3-\lambda_i\lambda_j}
{\lambda_i-\lambda_j}\right).
\end{eqnarray*}
Together with $b_1^2+b_2^2=1$, this gives a linear system of
three equations with unknowns $b_1^2$ and $b_2^2$. Since this
system is compatible by hypothesis, its rank is two and hence
the determinant of its augmented matrix is zero. This implies
\[
\frac{3c}{16(\lambda_1-\lambda_3)(\lambda_2-\lambda_3)}
(c-4\lambda_1\lambda_2+8\lambda_3(\lambda_1+\lambda_2)-12\lambda_3^2)=0.
\]
Solving the above system is only a matter of linear algebra.
After some calculations we get
$b_i^2=4(\lambda_j-2\lambda_3)(\lambda_i-\lambda_3)^2/(c(\lambda_i-\lambda_j))$
from where the result follows.
\end{proof}

We are now able to derive an important relation among $\lambda_1$, $\lambda_2$ and $\lambda_3$.

\begin{proposition}\label{thLambda12}
We have $c<0$. In this case, we get
\[
\lambda_i=\frac{1}{2}\left(3\lambda_3+(-1)^i\sqrt{-c-3\lambda_3^2}\right),
\quad (i,j\in\{1,2\},i\neq j).
\]
In particular, $\lambda_1<\lambda_3<\lambda_2$. Moreover,
$c+4\lambda_3^2<0$, or equivalently, $0\leq \lambda_3
<\sqrt{-c}/2$.
\end{proposition}

\begin{proof}
Let $i,j\in\{1,2\}$ with $i\neq j$. Using Lemma~\ref{thA}, the
constancy of $b_i$, and then Lemma~\ref{thNablaUiUjA}, we get
by Lemma~\ref{thGauss} applied to $U_i$ and $A$
\begin{eqnarray*}
0&=&
(\lambda_3-\lambda_i)(-c-4\lambda_i\lambda_3-2cb_j^2
+8\langle\nabla_{U_i}A,\nabla_A U_i\rangle)\\
&&-c b_i((-1)^ib_i\langle\nabla_A U_i,U_j\rangle
-2(-1)^j b_j\langle\nabla_{U_i}A,U_i\rangle
-2(-1)^i b_i\langle\nabla_{U_i}A,U_j\rangle)\\
&=&
\frac{c^2(\lambda_i-15\lambda_j+14\lambda_3)}
    {4(\lambda_3-\lambda_i)(\lambda_i-\lambda_j)}b_i^4
    +\frac{c^2(-10\lambda_i+3\lambda_j+7\lambda_3)}
    {2(\lambda_3-\lambda_i)(\lambda_i-\lambda_j)}b_i^2 b_j^2
    -\frac{11c\lambda_i(\lambda_3-\lambda_j)}{\lambda_i-\lambda_j}b_i^2\\
&&  -\frac{2c(\lambda_3-\lambda_i)(3\lambda_i-\lambda_j)}
    {\lambda_i-\lambda_j}b_j^2
    -\frac{(\lambda_3-\lambda_i)(c\lambda_i-c\lambda_j+8\lambda_i^2\lambda_j
    -4\lambda_3\lambda_i\lambda_j-4\lambda_3\lambda_i^2)}{\lambda_i-\lambda_j}.
\end{eqnarray*}
Now substituting $b_i^2$, $i\in\{1,2\}$, by the expressions
given in Proposition~\ref{thB}, after multiplying by $(\lambda_j-\lambda_i)/(\lambda_i-\lambda_3)$ and some long
calculations we get
\[
%\frac{\lambda_i-\lambda_3}{\lambda_i-\lambda_j}
72\lambda_3^3-48\lambda_i\lambda_3^2-108\lambda_j\lambda_3^2
+4\lambda_i^2\lambda_3+32\lambda_j^2\lambda_3
+72\lambda_i\lambda_j\lambda_3
-16\lambda_i\lambda_j^2-c\lambda_i
-8\lambda_i^2\lambda_j+c\lambda_j=0.
\]
Subtracting the above equation for $i=2$ from the one with
$i=1$ we get
$2(\lambda_1-\lambda_2)(c-4\lambda_1\lambda_2+14\lambda_3(\lambda_1+\lambda_2)
-30\lambda_3^2)=0$. Combining this with $c-4\lambda_1\lambda_2+8(\lambda_1+\lambda_2)\lambda_3-12\lambda_3^2=0$
(Proposition~\ref{thB}), we get
$6\lambda_3(\lambda_1+\lambda_2-3\lambda_3)=0$. If
$\lambda_3=0$, then the equation above gives
$c(\lambda_i-\lambda_j)+8\lambda_i^2\lambda_j+16\lambda_i\lambda_j^2=0$,
which combined with Proposition~\ref{thB} yields
$3c(\lambda_1+\lambda_2)=0$. This implies
$\lambda_1+\lambda_2-3\lambda_3=\lambda_1+\lambda_2=0$, so it suffices to deal with the case
$\lambda_1+\lambda_2-3\lambda_3=0$. In this situation we substitute
$\lambda_1$ by $-\lambda_2+3\lambda_3$ in the equation in
Proposition~\ref{thB}, thus obtaining
$c+4\lambda_2^2-12\lambda_2\lambda_3+12\lambda_3^2=0$. This is
a quadratic equation with unknown $\lambda_2$ and discriminant
$-c-3\lambda_3^2$. So that this discriminant is nonnegative we
already need $c<0$, proving the first claim of this
proposition. The solution to this equation is one of
\[
\frac{1}{2}\left(3\lambda_3\pm\sqrt{-c-3\lambda_3^2}\right).
\]
On the other hand, $\lambda_1$ is also one of the two values above.
Since $\lambda_1<\lambda_2$ by hypothesis, we get
$c+3\lambda_3^2<0$ and
$\lambda_i=\frac{1}{2}\left(3\lambda_3+(-1)^i\sqrt{-c-3\lambda_3^2}\right)$.

Finally, we show that $0\leq\lambda_3<\sqrt{-c}/2$. We already
know that $0\leq\lambda_3<\sqrt{-c/3}$. Substituting the above expression for $\lambda_i$, $i\in\{1,2\}$, in
Proposition~\ref{thB} we get
\[
b_i^2=-\frac{\left((-1)^i\lambda_3+\sqrt{-c-3\lambda_3^2}\right)^3}
{2c\sqrt{-c-3\lambda_3^2}},\quad i\in\{1,2\}.
\]
If $\sqrt{-c}/2\leq\lambda_3<\sqrt{-c/3}$, then
$-c-4\lambda_3^2\leq 0$, and hence
$-\lambda_3+\sqrt{-c-3\lambda_3^2}\leq 0$. This implies
$b_1^2\leq 0$, a contradiction. Therefore
$0\leq\lambda_3<\sqrt{-c}/2$ and the result follows.
\end{proof}

Proposition~\ref{thLambda12} already implies that there are no
hypersurfaces with constant principal curvatures in $\C
P^n(c)$, $n\geq 2$, whose Hopf vector field has $h=2$ nontrivial
projections onto the principal
curvature spaces. From now on we can assume $c<0$.

\begin{corollary}\label{thTlambdaReal}
The distribution $T_{\lambda_k}$ is totally real for all $k\in\{4,\dots,g\}$.
\end{corollary}

\begin{proof}
Let $k\in\{4,\dots,g\}$ and take unit vector fields $V_k,W_k\in\Gamma(T_{\lambda_k})$. Using the Weingarten equation, Lemma~\ref{thCodazzi}~(\ref{thCodazzi2Spaces}), Proposition~\ref{thB}, and $\lambda_1+\lambda_2-3\lambda_3=0$ (by Proposition~\ref{thLambda12}) we get
\begin{eqnarray*}
0&=&
V_k\langle W_k,J\xi\rangle\ =\
\langle\nabla_{V_k}W_k,b_1U_1+b_2U_2\rangle+\langle W_k,\bar\nabla_{V_k}J\xi\rangle\\
%&=&
%b_1\langle\nabla_{V_k}W_k,U_1\rangle+b_2\langle\nabla_{V_k}W_k,U_2\rangle
%-\lambda_k\langle JV_k,W_k\rangle\\
&=&
\left(\frac{cb_1^2}{4(\lambda_k-\lambda_1)}
+\frac{cb_2^2}{4(\lambda_k-\lambda_2)}-\lambda_k\right)\langle JV_k,W_k\rangle\ =\ \frac{(\lambda_3-\lambda_k)^3}{(\lambda_k-\lambda_1)(\lambda_k-\lambda_2)}
\langle JV_k,W_k\rangle.
\end{eqnarray*}
Since $\lambda_k\neq\lambda_3$, we get $\langle JV_k,W_k\rangle=0$. As $V_k$ and $W_k$ are arbitrary, the result follows.
\end{proof}

%%%%%%%%%%%%%%%%%%%%%%%%%% Subsection %%%%%%%%%%%%%%%%%%%%%%%%%%%%%%%
\subsection{A bound on the number of principal curvatures}\label{secBoundg}

In this section we show, using the Gauss equation and some
inequalities involving the principal curvatures, that the
number $g$ of distinct principal curvatures satisfies $g\in\{3,4\}$. This allows us to obtain further properties of the principal curvature spaces (see Proposition~\ref{thBoundg}). We start with the Gauss equation.

\begin{lemma}\label{thGaussInequalities}
Let us denote by $(\cdot)_i$, $i\in\{1,2\}$, the orthogonal
projection onto the distribution $T_{\lambda_i}\ominus\R U_i$,
and by $(\cdot)_k$, $k\in\{4,\dots,g\}$, the orthogonal
projection onto $T_{\lambda_k}$. By $\lVert\cdot\rVert$ we
denote the norm of a vector. Then we have:
\begin{enumerate}[{\rm (i)}]
\item\label{thGaussInequalities12} Let $i\in\{1,2\}$ and
    $W_i\in\Gamma(T_{\lambda_i}\ominus\R U_i)$ be a unit
    vector field. If $j\in\{1,2\}$ and $j\neq i$ then
    \[
    0=-(c+4\lambda_3\lambda_i)+8\frac{\lambda_i-\lambda_j}{\lambda_3-\lambda_j}
    \lVert(\nabla_A W_i)_j\rVert^2
    +8\sum_{k=4}^g\frac{\lambda_i-\lambda_k}{\lambda_3-\lambda_k}
    \lVert(\nabla_A W_i)_k\rVert^2.
    \]

\item\label{thGaussInequalitiesk} Let $k\in\{4,\dots,g\}$
    and $W_k\in\Gamma(T_{\lambda_k})$ be a unit vector field.
    Then
    \begin{eqnarray*}
    0&=&-(c+4\lambda_3\lambda_k)+8\frac{\lambda_k-\lambda_1}{\lambda_3-\lambda_1}
    \lVert(\nabla_A W_k)_1\rVert^2
    +8\frac{\lambda_k-\lambda_2}{\lambda_3-\lambda_2}
    \lVert(\nabla_A W_k)_2\rVert^2\\
    &&+8\sum_{l=4,\, l\neq k}^g\frac{\lambda_k-\lambda_l}{\lambda_3-\lambda_l}
    \lVert(\nabla_A W_k)_l\rVert^2.
    \end{eqnarray*}
\end{enumerate}
\end{lemma}

\begin{proof}
As usual, let $i,j\in\{1,2\}$ with $i\neq j$, and $k\in\{4,\dots,g\}$.

Let $W_i\in\Gamma(T_{\lambda_i})\ominus\R U_i$ be a unit
vector field. Applying Lemma~\ref{thGauss} to $W_i$ and $A$ we
get
\begin{equation}\label{eqGaussAWi}
-c-4\lambda_3\lambda_i+8\langle\nabla_{W_i}A,\nabla_A W_i\rangle=0.
\end{equation}
If $W_3\in\Gamma(T_{\lambda_3})$, we get from
Lemma~\ref{thCodazzi}~(\ref{thCodazzi2Spaces}) that
$\langle\nabla_A W_i,W_3\rangle=0$. This and
Lemma~\ref{thNablaUiUjA} yield $\nabla_A
W_i\in\Gamma((T_{\lambda_1}\ominus\R
U_1)\oplus(T_{\lambda_2}\ominus\R U_2)\oplus
T_{\lambda_4}\oplus\dots\oplus T_{\lambda_g})$. Similarly,
Lemma~\ref{thCodazzi}~(\ref{thCodazzi2Spaces}) implies
$\nabla_{W_i} A\in\Gamma(T_{\lambda_j}\oplus
(T_{\lambda_3}\ominus\R A)\oplus T_{\lambda_4}\oplus\dots\oplus
T_{\lambda_g})$. Hence $\langle\nabla_{W_i}A,\nabla_A
W_i\rangle=\langle\nabla_{W_i}A,(\nabla_A W_i)_j\rangle+
\sum_{k=4}^g\langle\nabla_{W_i}A,(\nabla_A W_i)_k\rangle$. For
each addend of this sum we apply
Lemma~\ref{thCodazzi}~(\ref{thCodazzi3Spaces}). Since
$\langle\bar{R}(W_i,A)(\nabla_A W_i)_l,\xi\rangle=0$ for all
$l\in\{j,4,\dots,g\}$ we get
\[
\langle\nabla_{W_i}A,\nabla_A W_i\rangle
=\frac{\lambda_i-\lambda_j}{\lambda_3-\lambda_j}
    \langle \nabla_A W_i,(\nabla_A W_i)_j\rangle
    +\sum_{k=4}^g\frac{\lambda_i-\lambda_k}{\lambda_3-\lambda_k}
    \langle\nabla_A W_i,(\nabla_A W_i)_k\rangle.
\]
Now, part~(\ref{thGaussInequalities12}) follows by substituting
the previous expression in \eqref{eqGaussAWi}.

Part (\ref{thGaussInequalitiesk}) follows in a similar way by
applying Lemma~\ref{thGauss} to $W_k$ and $A$.
\end{proof}

We will use the following technical lemma several times in what follows.

\begin{lemma}\label{thAuxiliaryLemma}
Assume $g\geq 4$ and let $k\in\{4,\dots,g\}$. Assume that one
of the following statements is true:
\begin{enumerate}[{\rm (i)}]
\item\label{thAuxiliaryLemmai} $\dim T_{\lambda_1}=\dim T_{\lambda_2}=1$, or

\item\label{thAuxiliaryLemmaii} $\dim T_{\lambda_1}=1$ and
    $\lambda_k<\lambda_2$, or

\item\label{thAuxiliaryLemmaiii} $\lambda_1<\lambda_k<\lambda_2$.
\end{enumerate}
Then, $c+4\lambda_3\lambda_k\geq 0$.
\end{lemma}

\begin{proof}
On the contrary, assume $c+4\lambda_3\lambda_k<0$. Let
$W_k\in\Gamma(T_{\lambda_k})$ be a (local) unit vector field.
When we apply
Lemma~\ref{thGaussInequalities}~(\ref{thGaussInequalitiesk}) to
$W_k$, any of the assumptions ensures that the first three
addends of the equation given in
Lemma~\ref{thGaussInequalities}~(\ref{thGaussInequalitiesk})
are nonnegative with the first one strictly positive. This
already implies $g>4$. In this case, it follows that there
exists $r\in\{4,\dots,g\}$, $r\neq k$, such that
$(\lambda_k-\lambda_r)/(\lambda_3-\lambda_r)<0$. We may choose
$\lambda_r$ to be the principal curvature that minimizes
$\lvert\lambda_3-\lambda_l\rvert$ among all $\lambda_l$,
$l\in\{4,\dots,g\}$, $l\neq k$, with
$(\lambda_k-\lambda_l)/(\lambda_3-\lambda_l)<0$. In particular
we have
\begin{equation}\label{eqLambdar}
\mbox{either }\quad\lambda_k<\lambda_r<\lambda_3\quad
\mbox{ or }\quad\lambda_3<\lambda_r<\lambda_k.
\end{equation}
It follows that $\lambda_r$ satisfies the same assumption as
$\lambda_k$: this is obvious for (\ref{thAuxiliaryLemmai}) and
a consequence of \eqref{eqLambdar} and
$\lambda_1<\lambda_3<\lambda_2$ for (\ref{thAuxiliaryLemmaii})
and (\ref{thAuxiliaryLemmaiii}). Using \eqref{eqLambdar},
$\lambda_3\geq 0$, $c+4\lambda_3^2<0$ (Proposition~\ref{thLambda12}), and
$c+4\lambda_3\lambda_k<0$, we also get
$c+4\lambda_3\lambda_r\leq
c+4\lambda_3\max\{\lambda_3,\lambda_k\}<0$. Thus we may apply
Lemma~\ref{thGaussInequalities}~(\ref{thGaussInequalitiesk}) to
a unit vector field $W_r\in\Gamma(T_{\lambda_r})$, from where it
follows, as before, that there exists $s\in\{4,\dots,g\}$,
$s\neq r$, such that
$(\lambda_r-\lambda_s)/(\lambda_3-\lambda_s)<0$. This implies
either $\lambda_r<\lambda_s<\lambda_3$ or
$\lambda_3<\lambda_s<\lambda_r$, and taking \eqref{eqLambdar}
into account we easily obtain
\begin{equation}\label{eqLambdars}
\mbox{either }\quad\lambda_k<\lambda_r<\lambda_s<\lambda_3\quad
\mbox{ or }\quad\lambda_3<\lambda_s<\lambda_r<\lambda_k.
\end{equation}
In both cases \eqref{eqLambdars} yields $s\neq k$,
$(\lambda_k-\lambda_s)/(\lambda_3-\lambda_s)<0$, and
$\lvert\lambda_3-\lambda_s\rvert<\lvert\lambda_3-\lambda_r\rvert$.
This contradicts the definition of $\lambda_r$. Therefore,
$c+4\lambda_3\lambda_k\geq 0$.
\end{proof}

From the previous lemma we easily derive the first important consequence.

\begin{proposition}\label{thDimT1}
We have $\dim T_{\lambda_1}=1$.
\end{proposition}

\begin{proof}
On the contrary, assume $\dim T_{\lambda_1}>1$ and let
$W_1\in\Gamma(T_{\lambda_1}\ominus\R U_1)$ be a (local) unit
vector field. Since $c+4\lambda_1\lambda_3\leq
c+4\lambda_3^2<0$ by Proposition~\ref{thLambda12}, from
Lemma~\ref{thGaussInequalities}~(\ref{thGaussInequalities12})
we deduce the existence of $k\in\{4,\dots,g\}$ such that
$(\lambda_1-\lambda_k)/(\lambda_3-\lambda_k)<0$. Since
$\lambda_1<\lambda_3$ we get
$\lambda_1<\lambda_k<\lambda_3<\lambda_2$ and hence
Lemma~\ref{thAuxiliaryLemma}~(\ref{thAuxiliaryLemmaiii}) yields
$c+4\lambda_3\lambda_k\geq 0$. This contradicts
$c+4\lambda_3\lambda_k\leq c+4\lambda_3^2<0$. Therefore $\dim
T_{\lambda_1}=1$.
\end{proof}

This is the most crucial step of the proof.

\begin{proposition}\label{thBoundg}
We have
\begin{enumerate}[{\rm (i)}]
\item\label{thBoundg34} $g\in\{3,4\}$.

\item\label{thBoundg3} If $g=3$ and $\dim T_{\lambda_2}>1$
    then $\lambda_1=0$, $\lambda_2=\frac{\sqrt{-3c}}{2}$,
    and $\lambda_3=\frac{\sqrt{-c}}{2\sqrt{3}}$.

\item\label{thBoundg4} If $g=4$ then $\dim
    T_{\lambda_2}=1$,  $0\neq\lambda_3\neq\frac{\sqrt{-c}}{2\sqrt{3}}$,
    and $\lambda_4=-\frac{c}{4\lambda_3}$.
\end{enumerate}
\end{proposition}

\begin{proof}
If $g=3$ and $\dim T_{\lambda_2}>1$, take
a (local) unit $W_2\in\Gamma(T_{\lambda_2}\ominus\R U_2)$ and apply
Lemma~\ref{thGaussInequalities}~(\ref{thGaussInequalities12}).
Note that the last two addends vanish since $\dim
T_{\lambda_1}=1$ and $g=3$. Then, $c+4\lambda_2\lambda_3=0$, and
from Proposition~\ref{thLambda12} we get $\lambda_1=0$,
$\lambda_2=\sqrt{-3c}/{2}$, and
$\lambda_3={\sqrt{-c}}/(2\sqrt{3})$. This implies
(\ref{thBoundg3}).

Assume $g\geq 4$. We first have $\lambda_3<\lambda_k$ for
all $k\in\{4,\dots,g\}$; otherwise, if
$\lambda_k<\lambda_3<\lambda_2$ we would get
$c+4\lambda_3\lambda_k\leq c+4\lambda_3^2<0$ contradicting
Lemma~\ref{thAuxiliaryLemma}~(\ref{thAuxiliaryLemmaii}) (by Proposition~\ref{thDimT1}).

We show that $\dim T_{\lambda_2}=1$. On the contrary,
assume $\dim T_{\lambda_2}>1$ and let
$W_2\in\Gamma(T_{\lambda_2}\ominus\R U_2)$ be a (local) unit
vector field. If $c+4\lambda_2\lambda_3<0$, then
Lemma~\ref{thGaussInequalities}~(\ref{thGaussInequalities12})
applied to $W_2$ (and taking Proposition~\ref{thDimT1} into account) implies that there exists $k\in\{4,\dots,g\}$
such that $(\lambda_2-\lambda_k)/(\lambda_3-\lambda_k)<0$. Then, $\lambda_3<\lambda_k<\lambda_2$, and thus
$c+4\lambda_3\lambda_k\leq c+4\lambda_3\lambda_2<0$, which
contradicts
Lemma~\ref{thAuxiliaryLemma}~(\ref{thAuxiliaryLemmaii}). Hence
we can assume from now on that $c+4\lambda_2\lambda_3\geq 0$.
This inequality does not hold if $\lambda_3=0$ so we already
get $\lambda_3>0$.

We claim that there exists $r\in\{4,\dots,g\}$ such that
$\lambda_2<\lambda_r$. If $c+4\lambda_2\lambda_3=0$, then the
assertion is true for all $k\geq 4$; otherwise, if
$\lambda_k<\lambda_2$, we would get $c+4\lambda_3\lambda_k<
c+4\lambda_3\lambda_2=0$, contradicting
Lemma~\ref{thAuxiliaryLemma}~(\ref{thAuxiliaryLemmaii}). Hence,
we have to prove our claim for the case
$c+4\lambda_2\lambda_3>0$. In this case we apply
Lemma~\ref{thGaussInequalities}~(\ref{thGaussInequalities12})
to $W_2$. Then, there exists  $r\in\{4,\dots,g\}$ such
that $(\lambda_2-\lambda_r)/(\lambda_3-\lambda_r)>0$. Since
$\lambda_3<\lambda_r$ this implies $\lambda_2<\lambda_r$ as
claimed.

In any case, there exists $r\in\{4,\dots,g\}$ such that
$\lambda_2<\lambda_r$. In fact, we may assume that $\lambda_r$
is the largest principal curvature. Now, we have
$c+4\lambda_3\lambda_r> c+4\lambda_3\lambda_2\geq 0$, and hence
Lemma~\ref{thGaussInequalities}~(\ref{thGaussInequalitiesk})
applied to a unit vector field $W_r\in\Gamma(T_{\lambda_r})$ implies
the existence of $l\in\{4,\dots,g\}$, $l\neq r$, such that
$(\lambda_r-\lambda_l)/(\lambda_3-\lambda_l)>0$. Since
$\lambda_3<\lambda_l$, we get $\lambda_r<\lambda_l$ which
contradicts the fact that $\lambda_r$ is the largest principal
curvature. Altogether this implies $\dim T_{\lambda_2}=1$.

From Lemma~\ref{thAuxiliaryLemma}~(\ref{thAuxiliaryLemmai}) we
obtain $c+4\lambda_3\lambda_k\geq 0$ for all $k\geq 4$. In
particular this implies $\lambda_3>0$. Assume that for some
$r\in\{4,\dots,g\}$ we have strict inequality
$c+4\lambda_3\lambda_r>0$ and let $\lambda_r$, $r\in\{4,\dots,g\}$, be the largest
principal curvature satisfying this condition. Applying
Lemma~\ref{thGaussInequalities}~(\ref{thGaussInequalitiesk})
once more to a unit $W_r\in\Gamma(T_{\lambda_r})$ (note that the
second addend now vanishes) yields the existence of
$l\in\{4,\dots,g\}$, $l\neq r$, such that
$(\lambda_r-\lambda_l)/(\lambda_3-\lambda_l)>0$. Since
$\lambda_3<\lambda_l$ we get $\lambda_r<\lambda_l$. Obviously,
$c+4\lambda_3\lambda_l>c+4\lambda_3\lambda_r>0$, which
contradicts the fact that $\lambda_r$ is the largest principal
curvature satisfying this condition.

As a consequence, $c+4\lambda_3\lambda_k=0$ for all $k\geq 4$.
Since $\lambda_3\neq 0$ and the principal curvatures are
different, this immediately implies $g=4$ and
$\lambda_4=-c/(4\lambda_3)$. Eventually, this also yields
$c+4\lambda_3\lambda_2\neq 0$ and thus, by
Proposition~\ref{thLambda12}, $\lambda_3\neq\sqrt{-c}/(2\sqrt{3})$ (otherwise the principal curvatures would not be different).
This concludes the proof of (\ref{thBoundg34}) and
(\ref{thBoundg4}).
\end{proof}

Part (\ref{thBoundg3}) of Proposition~\ref{thBoundg} had already
been obtained in \cite{BD05} by different methods. We have included a proof here as
it is almost effortless to do so.

%%%%%%%%%%%%%%%%%%%%%%%%%% Subsection %%%%%%%%%%%%%%%%%%%%%%%%%%%%%%%
\subsection{The eigenvalue structure of the shape
operator}\label{secEigenvalueStructure}

We summarize the results obtained so far:

\begin{theorem}\label{thSummary}
We have:
\begin{enumerate}[{\rm (a)}]
\item\label{thSummaryCPn} There are no real hypersurfaces with
    constant  principal curvatures in $\C P^n(c)$, $n\geq
    2$, whose Hopf vector field has $h=2$ nontrivial projections onto the principal curvature spaces.

\item\label{thSummaryCHn} Let $M$ be a connected real
    hypersurface with $g$ distinct constant principal curvatures
    $\lambda_1,\dots,\lambda_g$ in $\C H^n(c)$, $n\geq 2$,
    such that the number of nontrivial projections of its Hopf vector field $J\xi$ onto the
    principal curvature spaces of $M$ is $h=2$. Then,
    $g\in\{3,4\}$ and, with a suitable labeling of the
    principal curvatures and a suitable choice of the
    normal vector field $\xi$, we have:
    \begin{enumerate}[{\rm (i)}]
    \item The Hopf vector field can be written as
        $J\xi=b_1 U_1+b_2 U_2$, where
        $U_i\in\Gamma(T_{\lambda_i})$, $i\in\{1,2\}$,
        are unit vector fields, and $b_1$ and
        $b_2$ are positive constants satisfying
        \[
        b_i^2=\frac{4(\lambda_j-2\lambda_3)
            (\lambda_i-\lambda_3)^2}{c(\lambda_i-\lambda_j)},
            \quad (i,j\in\{1,2\},i\neq j).
        \]

    \item There exists a unit vector field
        $A\in\Gamma(T_{\lambda_3})$  such that
        \[
        JU_i    = (-1)^i b_j A-b_i\xi, \quad(i,j\in\{1,2\}, i\neq j),
        \quad\mbox{and}\quad
        JA      = b_2 U_1-b_1 U_1.
        \]

    \item We have $0\leq\lambda_3<\frac{1}{2}\sqrt{-c}$, and
        \[
        \lambda_i=\frac{1}{2}\left(3\lambda_3+(-1)^i\sqrt{-c-3\lambda_3^2}\right),
        \quad(i,j\in\{1,2\},i\neq j).
        \]

    \item $\dim T_{\lambda_1}=1$.

    \item\label{thSummaryCHng4} If $g=4$ then $\dim
        T_{\lambda_2}=1$.  We define $k=\dim
        T_{\lambda_4}+1$, and thus, $k\in\{2,\dots,n-1\}$.
        The distribution $T_{\lambda_4}$ is
        totally real with $JT_{\lambda_4}\subset
        T_{\lambda_3}\ominus\R A$,
        \[
        0\neq\lambda_3\neq\frac{\sqrt{-c}}{2\sqrt{3}},
        \quad\mbox{and}\quad \lambda_4=-\frac{c}{4\lambda_3}.
        \]

    \item\label{thSummaryCHng3} If $g=3$ there are two possibilities:
    \begin{enumerate}[{\rm (A)}]
    \item\label{thSummaryCHng31} $\dim
        T_{\lambda_2}=1$; in this  case we define
        $k=1$.

    \item\label{thSummaryCHng3k} $\dim
        T_{\lambda_2}>1$; in this  case we define
        $k=\dim T_{\lambda_2}\in\{2,\dots,n-1\}$
        and we have that $T_{\lambda_2}\ominus\R
        U_2$ is a real distribution with $J(T_{\lambda_2}\ominus\R
        U_2)\subset T_{\lambda_3}\ominus\R A$, and
        \[
        \lambda_1=0,\quad
        \lambda_2=\frac{\sqrt{-3c}}{2},\quad
        \lambda_3=\frac{\sqrt{-c}}{2\sqrt{3}}.
        \]
    \end{enumerate}
    \end{enumerate}
\end{enumerate}
\end{theorem}

\begin{remark}
Part~(\ref{thSummaryCPn}) of Theorem~\ref{thSummary} already
provides a proof for part~(\ref{thMainCPn}) of the Main
Theorem.

We know that $\R U_1\oplus\R U_2\oplus \R A\oplus\R\xi$ is a complex subbundle on $M$ by Lemma~\ref{thA}. Thus, in part (\ref{thSummaryCHng4}) of Theorem~\ref{thSummary}, the fact that $T_{\lambda_4}$ is real (Corollary~\ref{thTlambdaReal}) implies $JT_{\lambda_4}\subset T_{\lambda_3}\ominus\R A$ as claimed. Similarly, in Theorem~\ref{thSummary}~\ref{thSummaryCHng3k}, the assertion $J(T_{\lambda_2}\ominus\R U_2)\subset T_{\lambda_3}\ominus\R A$ follows from the fact that $T_{\lambda_2}$ is real by Lemma~\ref{thCodazzi}~(\ref{thCodazzi1Space}).

The definition of $k$ above might seem a bit artificial at the
moment, but it will be useful in the next section where we
conclude the proof of the Main Theorem ($k-1$ will be the dimension
of the kernel of the differential of the map $\Phi^r:M\to\C
H^n(c)$, $p\mapsto\exp_p(r\xi_p)$).

If we examine the proof of our theorem, so far we have actually
shown that for any point $p\in M$ there exists a neighborhood
of $p$ where the conclusion of Theorem~\ref{thSummary} is
satisfied. However, by the connectedness of $M$ and a
continuity argument, it can be easily shown that $M$ is
orientable and that the conclusion of Theorem~\ref{thSummary}
is satisfied globally.
\end{remark}

%%%%%%%%%%%%%%%%%%%%%%%%%% Subsection %%%%%%%%%%%%%%%%%%%%%%%%%%%%%%%
\subsection{Jacobi field theory and rigidity of focal
submanifolds}\label{secJacobi}

In this last section we finish the proof of part
(\ref{thMainCHn}) of the Main Theorem. Since we use standard
Jacobi field theory, we provide the reader just with the
fundamental details and skip the long calculations. According
to \cite{BD05} we just have to take care of the case $g=4$.
However, it is not much overload to deal with the two cases
simultaneously, so for the sake of completeness we will do so
in what follows.

Let $M$ be a real hypersurface of $\C H^n(c)$ in the conditions
of Theorem~\ref{thSummary}~(\ref{thSummaryCHn}). For $r\in\R$ we define the map
$\Phi^r:M\to\C H^n(c)$, $p\mapsto\exp_p(r\xi_p)$, where
$\exp_p$ is the Riemannian exponential map of $\C H^n(c)$ at
$p$. Then, $\Phi^r(M)$ is obtained by moving $M$ a distance $r$
along its normal direction. The singularities of $\Phi^r$ are
the focal points of $M$. We will find a particular distance $r$ for which $\Phi^r_*$ has
constant rank, where $\Phi^r_*$ denotes the differential of
$\Phi^r$. Then we will
apply Theorem~\ref{thRigidity} to $\Phi^r(M)$ for this choice of $r$ .
This way, $\Phi^r(M)$ will be an open part of the ruled minimal
Berndt-Br\"{u}ck submanifold $W^{2n-k}$, $k\in\{1,\dots,n-1\}$, and hence $M$
will be an open part of a tube around this ruled minimal
submanifold  $W^{2n-k}$. (If $k=1$ then $M$ will be an
equidistant hypersurface to the ruled minimal hypersurface
$W^{2n-1}$ at distance $r$.)

Let $p\in M$ and denote by $\gamma_p$ the geodesic determined
by the initial conditions $\gamma_p(0)=p$ and
$\dot\gamma_p(0)=\xi_p$. For any $v\in T_p M$ let $B_v$ be the
parallel vector field along the geodesic $\gamma_p$ such that
$B_v(0)=v$, and let $\zeta_v$ be the Jacobi field along
$\gamma_p$ with initial conditions $\zeta_v(0)=v$ and
$\zeta'(0)=-S_p v$. Here $'$ denotes covariant derivative along
$\gamma_p$. Since $\zeta_v$ is a solution to the differential
equation
$4\zeta_v''-c\zeta_v-3c\langle\zeta_v,J\dot\gamma_p\rangle
J\dot\gamma_p=0$, if $v\in T_{\lambda_i}(p)$ then
\[
\zeta_v(t)=f_i(t)B_v(t)+\langle v,J\xi\rangle g_i(t)J\dot\gamma_p(t),
\]
where
\begin{eqnarray*}
f_i(t)&=&
\cosh\left(\frac{t\sqrt{-c}}{2}\right)-\frac{2\lambda_i}{\sqrt{-c}}
\sinh\left(\frac{t\sqrt{-c}}{2}\right),\\
g_i(t)&=&
\left(\cosh\left(\frac{t\sqrt{-c}}{2}\right)-1\right)
\left(1+2\cosh\left(\frac{t\sqrt{-c}}{2}\right)
-\frac{2\lambda_i}{\sqrt{-c}}\sinh\left(\frac{t\sqrt{-c}}{2}\right)\right).
\end{eqnarray*}
We also define the smooth vector field $\eta^r$ along $\Phi^r$
by $\eta^r_p=\dot\gamma_p(r)$. It is known that
$\zeta_v(r)=\Phi^r_*v$ and
$\zeta_v'(r)=\bar\nabla_{\Phi^r_*v}\eta^r$.

We now determine the value of $r$. Since $0\leq
\lambda_3<\sqrt{-c}/2$  we can find a real number $r\geq 0$
such that
\[
\lambda_3=\frac{\sqrt{-c}}{2}\tanh\left(\frac{r\sqrt{-c}}{2}\right).
\]
%If $r=0$ then $\Phi^0$ is just the inclusion map.

Let $p\in M$. We define $u_i=(U_i)_p$, $i\in\{1,2\}$. Let
$v_2\in T_{\lambda_2}(p)\ominus\R u_2$ and $v_k\in
T_{\lambda_k}(p)$ for $3\leq k\leq g$ (whenever these spaces
are nontrivial). The explicit solution to the Jacobi equation
above implies
\[
(\Phi^r_* u_1, \Phi^r_* u_2)=(B_{u_1}(r), B_{u_2}(r)) D(r),
\]
\[
\Phi^r_* v_2 = 0, \qquad
\Phi^r_* v_3 = \sech\left(\frac{t\sqrt{-c}}{2}\right)B_{v_3}(r),\qquad
\Phi^r_* v_4 = 0,
\]
where
\[
D(t)=\left(\begin{array}{cc}
f_1(t)+b_1^2 g_1(t) &   b_1b_2 g_2(t)\\
b_1b_2 g_1(t)   &   f_2(t)+b_2^2 g_2(t)
\end{array}\right).
\]
Since $\det(D(r))=\sech^3\left({r\sqrt{-c}}/{2}\right)$ we
conclude that $\Phi^r_*$ has constant rank $2n-k$ (see
Theorem~\ref{thSummary}~(\ref{thSummaryCHng4})-(\ref{thSummaryCHng3})
for the definition of $k$). Then, for each point $p\in M$ there
exists an open neighborhood $\mathcal{V}$ of $p$ such that
$\mathcal{W}=\Phi^r(\mathcal{V})$ is an embedded submanifold of
$\C H^n(c)$ and $\Phi^r:\mathcal{V}\to\mathcal{W}$ is a
submersion. (If $k=1$, then $\Phi^r$ is actually a local
diffeomorphism.)

Let $q=\Phi^r(p)\in\mathcal{W}$. The expression above for
$\Phi^r_*$ shows that the tangent space $T_q\mathcal{W}$ of
$\mathcal{W}$ at $q$ is obtained by parallel translation of $\R
u_1\oplus\R u_2\oplus T_{\lambda_3}(p)$ along the geodesic
$\gamma_p$ from $p=\gamma_p(0)$ to $q=\gamma_p(r)$. Therefore,
the normal space $\nu_q\mathcal{W}$ of $\mathcal{W}$ at $q$ is
obtained by parallel translation of $(\ker \Phi^r_{*p})\oplus\R\xi_p$
along $\gamma_p$ from $p=\gamma_p(0)$ to $q=\gamma_p(r)$. The latter
is $(T_{\lambda_2}\ominus\R u_2)\oplus\R\xi_p$ if $g=3$
(see Theorem~\ref{thSummary}~(\ref{thSummaryCHng3})), or
$T_{\lambda_4}(p)\oplus\R\xi_p$ if $g=4$ (see
Theorem~\ref{thSummary}~(\ref{thSummaryCHng4})). In any case,
by Theorem~\ref{thSummary}~(\ref{thSummaryCHng4})-(\ref{thSummaryCHng3})
it follows that $\mathcal{W}$ has totally real normal bundle of
rank $k$.

We have that $\eta^r_p=B_{\xi_p}(r)$ is a unit normal vector
of $\mathcal{W}$ at $q$. If $S^r$ denotes the shape operator of
$\mathcal{W}$, then it is known that
$S^r_{\eta^r_p}\Phi^r_*v=-(\zeta_v'(r))^\top$, where
$(\cdot)^\top$ denotes orthogonal projection onto the tangent
space of $\mathcal{W}$. Using the explicit expression for
$\zeta_v$ above, we get
\[
(S^r_{\eta^r_p}B_{u_1}(r), S^r_{\eta^r_p}B_{u_2}(r))
=(B_{u_1}(r), B_{u_2}(r)) C(r),
\quad\mbox{and}\quad
S^r_{\eta^r_p}B_{v_3}(r)=0\mbox{  for all } v_3\in T_{\lambda_3}(p),
\]
where $C(r)=-D'(r)D(r)^{-1}$. A lengthy and tedious calculation shows that
\[
C(r)=\frac{\sqrt{-c}}{2}\left(\begin{array}{cc}
-2b_1b_2        &   b_1^2-b_2^2\\
b_1^2-b_2^2     &   2b_1b_2
\end{array}\right).
\]
Since $J\eta^r_p=B_{J\xi_p}(r)=b_1 B_{u_1}(r)+b_2 B_{u_2}(r)$,
and $B_{J A_p}(r)=b_2 B_{u_1}(r)-b_1 B_{u_2}(r)$, the above
expression for $C(r)$ implies
\[
S^r_{\eta^r_p}B_{J A_p}(r)=-\frac{\sqrt{-c}}{2}J\eta^r_p,\qquad
S^r_{\eta^r_p}J\eta^r_p=-\frac{\sqrt{-c}}{2}B_{J A_p}(r),
\]
and $S^r_{\eta^r_p}$ vanishes on the orthogonal complement of
$\R J\eta^r_p\oplus\R B_{J A_p}(r)$ in $T_{q}\mathcal{W}$.

We have that $J(\nu_q \mathcal{W}\ominus\R\eta^r_p)$ is
contained in the parallel translation along $\gamma_p$ of $T_{\lambda_3}(p)$.
This follows from
Theorem~\ref{thSummary}~(\ref{thSummaryCHng4})-(\ref{thSummaryCHng3})
and the fact that $\nu_q\mathcal{W}\ominus\R\eta^r_p$ is the
parallel translation along $\gamma_p$ from $\gamma_p(0)=p$ to
$\gamma_p(r)=q$ of $T_{\lambda_2}(p)\ominus\R u_2$ if $g=3$, and
of $T_{\lambda_4}(p)$ if $g=4$. The linearity of
$S^r_{\eta^p_r}$ implies
\begin{equation}\label{eqSJtildexi}
S^r_{\eta^r_p}J\tilde\eta=-\frac{\sqrt{-c}}{2}\langle\eta^r_p,
\tilde\eta\rangle B_{JA_p}(r),
\mbox{ for all } \tilde\eta\in\nu_q\mathcal{W}.
\end{equation}
It follows from the Gauss formula and $\bar\nabla J=0$ that
$S^r_{\tilde\eta}J\eta^r_p=S^r_{\eta^r_p}J\tilde\eta$, and hence,
$S^r_{\tilde\eta}J\eta^r_p=0$ for all
$\tilde\eta\in\nu_q\mathcal{W}\ominus\R\eta^r_p$. Let $\alpha$
be a curve in $(\Phi^r)^{-1}(\{q\})\cap\mathcal{V}$ with
$\alpha(0)=p$. Since $\eta^r_p$ and
$\eta^r_{\alpha(t)}-\langle\eta^r_{\alpha(t)},\eta^r_p\rangle\eta^r_p$
are perpendicular, $S^r_{\tilde\eta}J\eta^r_p=0$, and the
linearity of $\eta\mapsto S^r_\eta$ imply
\[
0=S^r_{\eta^r_{\alpha(t)}-\langle\eta^r_{\alpha(t)},
\eta^r_p\rangle\eta^r_p}J\eta^r_p=
S^r_{\eta^r_{\alpha(t)}}J\eta^r_p
+\frac{\sqrt{-c}}{2}\langle\eta^r_{\alpha(t)},\eta^r_p\rangle B_{JA_p}(r),
\]
which together with \eqref{eqSJtildexi} (with $\alpha(t)$
instead of $p$) yields
\[
-\frac{\sqrt{-c}}{2}\langle\eta^r_{\alpha(t)},\eta^r_p\rangle B_{JA_p}(r)=
S^r_{\eta^r_{\alpha(t)}}J\eta^r_p
=-\frac{\sqrt{-c}}{2}\langle\eta^r_{\alpha(t)},\eta^r_p\rangle
B_{JA_{\alpha(t)}}(r).
\]
Since $\alpha$ is arbitrary we get that the map
$\tilde{p}\mapsto B_{JA_{\tilde{p}}}(r)$ is constant in the
connected component $\mathcal{V}_0$ of
$(\Phi^r)^{-1}(\{q\})\cap\mathcal{V}$ containing $p$. Thus it
makes sense to define the unit vector
$z=-B_{JA_{\tilde{p}}}(r)\in T_q\mathcal{W}$ for any
$\tilde{p}\in\mathcal{V}_0$.

We may consider $\eta^r$  as a map from $\mathcal{V}_0$ to the
unit sphere of $\nu_q\mathcal{W}$. The tangent space of
$\mathcal{V}_0$ at $p$ is given by the kernel of $\Phi^r_{*p}$. If
$v\in\ker\Phi^r_{*p}$, then $\eta^r_{*p}v=\zeta_v'(r)$. If
$g=3$, then $v\in\ker\Phi^r_{*p}=T_{\lambda_2}(p)\ominus\R u_2$,
and $\eta^r_{*p}=-\sqrt{-c/2}\,B_{v}(r)$. If $g=4$, then
$v\in\ker\Phi^r_{*p}=T_{\lambda_4}(p)$, and
$\eta^r_{*p}v=-\csch(r\sqrt{-c}/2)B_{v}(r)$. In any case, we get
that $\eta^r$ is a local diffeomorphism from $\mathcal{V}_0$
into the unit sphere of $\nu_q\mathcal{W}$ (note that this is
trivial if $g=3$ and $k=1$). Hence, $\eta^r(\mathcal{V}_0)$ is
an open subset of the unit sphere of $\nu_q\mathcal{W}$. But
since $\eta\mapsto S^r_\eta$ depends analytically on $\eta$ we
conclude
\[
S^r_\eta J\eta=\frac{\sqrt{-c}}{2}z,\qquad
S^r_\eta z=\frac{\sqrt{-c}}{2}J\eta,\qquad
S^r_\eta v=0,
\]
for all unit $\eta\in\nu_q\mathcal{W}$, and $v\in
T_q\mathcal{W}\ominus(\R J\eta\oplus\R z)$. Therefore, the
second fundamental form $\II^r$ of $\mathcal{W}$ at $q$ is
given by the trivial symmetric bilinear extension of
$\II^r(z,J\eta)=(\sqrt{-c}/2)\eta$ for all
$\eta\in\nu_q\mathcal{W}$. By construction, $z$ depends smoothly
on the point $q\in\mathcal{W}$ and hence gives rise to a vector
field $Z$ which is tangent to the maximal holomorphic
distribution of $\mathcal{W}$. The relation
$S^r_{\eta}J\eta=(\sqrt{-c}/2)Z$ ensures that $Z$ can actually
be defined on $\Phi^r(M)$, and hence, the second fundamental form
of $\Phi^r(M)$ is given by the trivial symmetric bilinear extension of
$\II^r(Z,J\eta)=(\sqrt{-c}/2)\eta$ for all
$\eta\in\Gamma(\nu\,\Phi^r(M))$. Since $\Phi^r(M)$ has totally
real normal bundle of rank $k$ we conclude from
Theorem~\ref{thRigidity}, and the remark that follows, that $\Phi^r(M)$ is holomorphically
congruent to an open part of the ruled minimal Berndt-Br\"{u}ck submanifold
$W^{2n-k}$. This readily implies that $M$ is an open part of a
tube (an equidistant hypersurface if $g=3$ and $k=1$) of radius
$r$ around the ruled minimal Berndt-Br\"{u}ck submanifold $W^{2n-k}$.

Finally, let us point out that if $g=3$ and $\lambda_3=0$, then
$r=0$ and $M$ is an open part of the ruled minimal hypersurface
$W^{2n-1}$. Also, if $g=3$ and $k>1$ then
$\lambda_3=\sqrt{-c}/(2\sqrt{3})$ according to
Theorem~\ref{thSummary}~\ref{thSummaryCHng3k}, and hence
$r=(1/\sqrt{-c})\log(2+\sqrt{3})$. The tube around the ruled
minimal submanifold  $W^{2n-k}$, $k>1$, of radius
$r=(1/\sqrt{-c})\log(2+\sqrt{3})$ has $g=3$ principal
curvatures whereas if $r\neq(1/\sqrt{-c})\log(2+\sqrt{3})$ the
tube of radius $r$ around the ruled minimal submanifold
$W^{2n-k}$, $k>1$, has $g=4$ principal curvatures.

This finishes the proof of the Main Theorem.

%%%%%%%%%%%%%%%%%%%%%%%%%% Bibliography %%%%%%%%%%%%%%%%%%%%%%%%%%%

\end{document}